\documentclass[reqno,10pt]{amsart}  

\usepackage{amssymb,amsmath}
\usepackage{amsthm}  
\usepackage{latexsym}
\usepackage{amscd}
\usepackage{color}  
\usepackage{fancybox, framed}  
\usepackage[normalem]{ulem}  
\usepackage{eufrak}  

\let\cal=\mathcal

\def\bbAbip{\bbA_{\text{{\rm bip}}}}
\def\bbAappbip{\bbA_{\approx\text{{\rm bip}}}}

\def\bbAbipelem{\bbA_{\text{{\rm elem-bip}}}}
\def\bbAbipreg{\bbA_{\text{{\rm reg-bip}}}}

\def\cXH{{\cal X}_H}  
\def\cYH{{\cal Y}_H}

\def\Xapp{X^{\approx}}  
\def\Yapp{Y^{\approx}}  

\def\XappH{X^{\approx}_H}  
\def\YappH{Y^{\approx}_H}  

\def\XedH{X^{\rm ed}_H}  
\def\YedH{Y^{\rm ed}_H}  

\def\Xed{X^{\rm ed}}  
\def\Yed{Y^{\rm ed}}

\def\dcup{\, \dot\cup \,}

\def\epsilon{\varepsilon}

\def\eps{\varepsilon}

\def\cB{{\cal B}}

\def\cE{{\cal E}}
\def\cF{{\cal F}}

\def\cJ{{\cal J}}

\def\cP{{\cal P}}
\def\cQ{{\cal Q}}
\def\cR{{\cal R}}
\def\cS{{\cal S}}
\def\cT{{\cal T}}

\def\cX{{\cal X}}
\def\cY{{\cal Y}}

\def\dcup{\, \dot{\cup} \,}  

\def\bbB{\mathbb{B}}  
\def\bbE{\mathbb{E}}  
\def\bbH{\mathbb{H}}  
\def\bbP{\mathbb{P}}  

\def\bC{\boldsymbol{C}}

\def\bbA{\mathbb{A}}

\usepackage{amsthm,amsmath,amssymb}

\title{On two-coloring bipartite uniform hypergraphs}    

\author{Boyoon Lee}
\address{Department of Mathematics and Statistics, 
University of South Florida, 
Tampa, FL 33620, USA.}
\email{jihyelee@usf.edu}

\author{Theodore Molla}\thanks{The second author was partially supported
by NSF Grants DMS~1800761 and DMS~2154313.}  
\address{Department of Mathematics and Statistics, 
University of South Florida, 
Tampa, FL 33620, USA.}
\email{molla@usf.edu}

\author{Brendan Nagle}\thanks{The third author was partially supported by 
NSF Grant DMS~1700280.}  
\address{Department of Mathematics and Statistics, 
University of South Florida, 
Tampa, FL 33620, USA.}
\email{bnagle@usf.edu}

\usepackage{amsthm,amsmath,amssymb}
\usepackage{latexsym}
\usepackage{amsmath}
\usepackage{amssymb}
\usepackage{amsthm}
\usepackage{amsfonts}

\theoremstyle{plain}  

\newtheorem{theorem}{Theorem}[section]

\newtheorem{lemma}[theorem]{Lemma}
\newtheorem{fact}[theorem]{Fact}

\theoremstyle{definition}
\newtheorem{definition}[theorem]{Definition}
\newtheorem{remark}[theorem]{Remark}

\newtheorem{cor}[theorem]{Corollary}  
\newtheorem{innercustomthm}{Theorem}
\newenvironment{customthm}[1]
  {\renewcommand\theinnercustomthm{#1}\innercustomthm}
  {\endinnercustomthm}

\newtheorem*{remark*}{Remark}



\begin{document}

\maketitle

\begin{abstract}
Of a given bipartite graph $G = (V, E)$, 
it is elementary to construct a bipartition in time $O(|V| + |E|)$.  
For a given $k$-graph 
$H = H^{(k)}$ with $k \geq 3$ fixed, Lov\'asz proved that deciding whether $H$ is bipartite
is NP-complete.  
Let $\cB_n$ denote the collection of all $[n]$-vertex bipartite $k$-graphs.  
We construct, of a given $H \in \cB_n$, 
a bipartition in time averaging $O(n^k)$ over the class $\cB_n$.  
We provide two proofs of our result.  
When $k = 3$, this result expedites one of Person and Schacht.  
\end{abstract}

\section{Introduction}  
Fix an integer $k \geq 2$.
A {\it $k$-uniform hypergraph} $H = H^{(k)}$, or {\it $k$-graph} for short,  
is a pair $(V_H, E_H)$ where $E_H \subseteq \binom{V_H}{k}$  
is a family of $k$-element subsets from $V_H$.  
We say $H$ 
is {\it bipartite} when there exists a vertex partition $V_H = X_H \, \dcup \, Y_H$, 
which we call a {\it bipartition of $H$}  
(or a {\it 2-coloring of $H$}), 
where both $X_H$ and $Y_H$ meet every $e \in E_H$.  
When $k = 2$, 
it is elementary to find in time $O(|V_H| + |E_H|)$ either a bipartition
of $H$ or an odd cycle therein.  For $k \geq 3$, Lov\'asz~\cite{L} proved
that deciding whether $H$ is bipartite is NP-complete.  
For $k \geq 4$, 
Guruswami, H\aa stad, and Sudan~\cite{GHS} proved that it is NP-hard to properly color
a given bipartite $k$-graph 
$H$ with constantly many colors.

In this paper, 
we consider 2-coloring bipartite $k$-graphs in terms of average running time.  
For that, let 
$\cB_n = \cB_n^{(k)}$ be the class of all 
$n$-vertex  
bipartite $k$-graphs $H$ with the fixed and labelled vertex 
set $V_H = [n] = \{1, \dots, n\}$.   
Throughout this paper, 
we take $n \in \mathbb{N}$ sufficiently large wherever needed.  

\begin{customthm}{A}
\label{thm:main}  
{\it  
There exists an algorithm $\bbA_{{\rm bip}}$
which constructs, for each fixed integer $k \geq 3$ 
and 
of a given $k$-graph $H \in \cB_n$, 
a bipartition 
$[n] = X_H \dcup Y_H$ 
of $H$ 
in time averaging $O(n^k)$ over the class $\cB_n$.}  
\end{customthm}

\noindent  In this paper, we will provide two proofs of Theorem~\ref{thm:main}.  
We postpone discussion
on these proofs 
for a moment to establish further context
for Theorem~\ref{thm:main}.  

Theorem~\ref{thm:main} follows a tradition of related results. 
Dyer and Frieze~\cite{DF} established an algorithm properly $\ell$-coloring 
a given $\ell$-colorable $[n]$-vertex graph $G$ in average time $O(n^2)$
(see also~\cite{PrSt, T, W}).  
Person and Schacht~\cite{PSalm, PSalg} 
established a hypergraph analogue thereof 
regarding the {\it Fano plane $F$}, which is 
the unique $(7,3,1)$-design.  
An elementary property of 
the Fano plane $F$ is that it is not bipartite.  
In what follows, we let 
$\cF_n = \cF_n^{(3)}$
be the class of all 
$[n]$-vertex 3-graphs $H$ forbidding the Fano plane $F$ as a subhypergraph.    
Then 
$\cB_n^{(3)} \subseteq \cF_n$, which Person and Schacht~\cite{PSalm, PSalg} showed is near equality.  

\begin{theorem}[Person and Schacht~\cite{PSalm, PSalg}]  
\label{thm:PSFano}  
The inequality 
$\big|\cF_n \setminus \cB_n^{(3)}\big| \leq \big|\cB_n^{(3)}\big| / 2^{\Omega (n^2)}$ holds.  
Moreover,  
there exists an algorithm $\bbA_{\rm Fano}$ 
which constructs, of a given $H \in \cF_n$, 
a smallest partition 
$[n] = X_{H,1} \dcup \dots \dcup X_{H,r}$ 
into $r = r_H = \chi(H)$ 
$($the weak chromatic number of $H)$
independent sets of $H$,     
in time averaging $O(n^5 \log^2 n)$ over $\cF_n$.    
\end{theorem}  

Person and Schacht~\cite{PSalm, PSalg} note the following corollary of Theorem~\ref{thm:PSFano} and precursor 
of Theorem~\ref{thm:main}.  

\begin{cor}[Person and Schacht~\cite{PSalm, PSalg}]
\label{cor:PS}  
The algorithm $\bbA_{{\rm Fano}}$ constructs, of a given $H \in \cB_n^{(3)}$, 
a bipartition $[n] = X_H \dcup Y_H$ in time averaging $O(n^5 \log^2 n)$ over 
the class 
$\cB_n^{(3)}$.  
\end{cor}  

\noindent  Corollary~\ref{cor:PS} 
inspires Theorem~\ref{thm:main}, 
but the main ideas of the proofs differ.  
Corollary~\ref{cor:PS} 
depends on constructive hypergraph regularity and counting lemmas
from~\cite{CR, KNRS}, which  
are needed to prove the stronger Theorem~\ref{thm:PSFano}.    
We give two proofs of Theorem~\ref{thm:main} avoiding these tools.
The first is elementary and based 
on the simple tools 
of convexity and the Chernoff inequality $\big($cf.~\cite{JLR}$\big)$.    
It shows that
all but $|\cB_n| / 2^{\Omega(n^{k-1})}$ many $H \in \cB_n$ 
admit a unique bipartition 
$[n] = X_H \dcup Y_H$ 
constructible in maximum time $O(n^k)$.     
The second is based on 
the constructive Szemer\'edi Regularity Lemma~\cite{KRT, 
Szem1, Szem} and 
gives a structural nuance of possible interest.  It constructs, for arbitrary $k \geq 3$
but in time just $O(n^2)$, a partition 
$[n] = X_H^{\approx} \dcup Y_H^{\approx}$ 
for $H \in \cB_n$ which, for all but 
$|\cB_n| / 2^{\Omega(n^2)}$ many $H \in \cB_n$, closely approximates  
the unique bipartition of $H$ in the sense that 
$\big|X_H \, \triangle \, X_H^{\approx}\big| = 
\big|Y_H \, \triangle \,  Y_H^{\approx}\big| = o(n)$.

\subsection*{Itinerary of paper}      
Section~\ref{sec:shortproof}  
presents the elementary proof  
of Theorem~\ref{thm:main}, and denotes by 
$\bbAbipelem$ 
the version of $\bbAbip$ corresponding to this proof.
Section~\ref{sec:bbAbipreg}
presents the 
regularity-based 
version of $\bbAbip$, which it denotes by $\bbAbipreg$.  
Section~\ref{sec:mainproofsep} 
proves 
Theorem~\ref{thm:main}
from $\bbAbipreg$.
Section~\ref{sec:firstlemma} proves an important lemma needed along the way.    
Section~\ref{sec:Fano}  
closes with some concluding remarks on 
Theorem~\ref{thm:PSFano}.    
The Appendix proves a few standard but technical facts needed in this paper.

\subsubsection*{Acknowledgment}  
The authors would like to thank their highly diligent Referees, whose excellent suggestions led to a much improved
presentation of this paper.  The authors would also like to thank Marcus Michelen 
for 
mentioning 
an idea to us 
that would ultimately inspire
Section~2.     
We also thank him for his selfless encouragement to pursue this work.

\section{An Elementary Proof of Theorem~\ref{thm:main}}  
\label{sec:shortproof}  
In this section, we present the elementary proof of Theorem~\ref{thm:main} for a fixed
integer $k \geq 3$.   
To do so, we first present its algorithm, 
denoted here by 
$\bbAbipelem$.   
We begin our work by establishing notation.

\subsection*{Notation}  
Fix 
an $[n]$-vertex $k$-graph 
$H$  
and distinct vertices $u \neq v \in [n]$.  
We 
define 
\begin{equation}
\label{eqn:basicneigh}  
N_H(u) = \big\{J \subseteq [n] \setminus \{u\} : \, |J| = k - 1 
\text{ and } \{u\} \dcup J \in E_H \big\}  
\end{equation}  
to be the 
{\it $H$-neighborhood of $u$}.
We respectively define
\begin{equation}
\label{eqn:jointdeg}
N_H(u) \cap 
N_H(v) 
\qquad \text{and} \qquad 
\deg_H(u, v) 
= \big|
N_H(u) \cap 
N_H(v)\big| 
\end{equation}  
to be 
the {\it joint $H$-neighborhood} and 
{\it joint $H$-degree of $u$ and $v$}.

\begin{remark*}
We pause for a word on~\eqref{eqn:jointdeg}.  Note that 
joint $H$-degrees $\deg_H(u, v)$ 
count $(k-1)$-tuples $J$ 
for which both $\{u\} \dcup J \in E_H$ and $\{v\} \dcup J \in E_H$ hold.   
In contrast, 
{\it $H$-codegrees} $\deg_H\big(\{u, v\}\big)$  
(which are quite common in the literature)  
count $(k-2)$-tuples $I$
for which $\{u, v\} \dcup I \in E_H$.  
\end{remark*}

Joint $H$-degrees~\eqref{eqn:jointdeg}
are the central consideration of $\bbAbipelem$.  
Before we present those details, 
we motivate and outline them using 
a standard random hypergraph model which will be relevant for Theorem~\ref{thm:main}.  
For that, fix a partition\footnote{In this paper, $o(1) \to 0^+$ as $n \to \infty$.}
$[n] = X \dcup Y$ satisfying 
$|X|, |Y| = \big(1 \pm o(1) \big) n / 2$ and 
consider 
the {\it $[n]$-vertex binomial random $k$-graph 
$\mathbb{H}_{X, Y} = \mathbb{H}^{(k)}[X, Y]$}
which 
independently includes each $K \in \binom{[n]}{k} \setminus \big(\binom{X}{k} \dcup \binom{Y}{k} \big)$ 
as an edge with probability
$1/2$.  
We proceed with some standard but motivating details.  

\subsection*{Expected joint $\boldsymbol{\mathbb{H}_{X,Y}}$-degrees
and $\boldsymbol{\mathbb{A}_{\rm elem-bip}}$}  
Fix $u \neq v \in [n]$.  
Clearly, 
\begin{equation}
\label{eqn:expjointpure}  
\mathbb{E} 
\big[ \deg_{\bbH_{X,Y}}(u, v) \big]
 =  
\left\{
\begin{array}{ll}
\tfrac{1}{4} \big(\tbinom{n-2}{k-1} - \tbinom{|X|-2}{k-1} \big)
& 
\text{when $u, v \in X$,} \bigskip \\
\tfrac{1}{4} \big(\tbinom{n-2}{k-1} - \tbinom{|Y|-2}{k-1} \big)
& 
\text{when $u, v \in Y$,} \bigskip \\
\tfrac{1}{4} \big(\tbinom{n-2}{k-1} 
- \tbinom{|X|-1}{k-1} 
- \tbinom{|Y|-1}{k-1} \big)
& \text{otherwise.}
\end{array}
\right.  
\end{equation}  
From $|X|, |Y| = \big(1 \pm o(1) \big) n / 2$ follows
\begin{equation}
\label{eqn:expjoint}  
\mathbb{E} 
\big[ \deg_{\bbH_{X,Y}}(u, v) \big]
 =  
\left\{
\begin{array}{ll}  
\tfrac{1}{4} 
\tbinom{n}{k-1} \big(1 - \tfrac{1}{2^{k-1}} \pm o(1) \big)    
& 
\text{when $u, v\in X$ or $u, v \in Y$,} \bigskip \\
\tfrac{1}{4} \tbinom{n}{k-1} \big(1 - \tfrac{2}{2^{k-1}} \pm o(1) \big)    
& \text{otherwise.}
\end{array}
\right.
\end{equation}

Fix now $H \in \cB_n$.
The algorithm $\bbAbipelem$ will first construct a partition $[n] = X^H \dcup Y^H$    
which it intends (however incorrectly) to be a bipartition of $H$.  
It will decide to 
place both $u, v \in X^H$ or $u, v \in Y^H$ 
if, and only if, 
$\deg_H(u, v)$ 
is closer 
to 
$(1/4) \tbinom{n}{k-1} \big(1 - 2^{1-k} \big)$
$\big($cf.~\eqref{eqn:expjoint}$\big)$   
than it is to 
$(1/4) \tbinom{n}{k-1} \big(1 - 2^{2-k} \big)$
$\big($cf.~\eqref{eqn:expjoint}$\big)$, which happens 
if, and only if, 
$\deg_H(u, v)$ 
is at least their average:  
\begin{equation}
\label{eqn:1.17.2025.1:14p}  
\deg_H(u, v) 
\geq 
\tfrac{1}{4} \tbinom{n}{k-1} \big(1 - \tfrac{3}{2^k}\big).  
\end{equation}  
In particular, 
$\bbAbipelem$ will iteratively   
build $[n] = X^H \dcup Y^H$
by the rule that if, w.l.o.g., 
$\bbAbipelem$ last placed $u \in X^H$, then it next places $v \in X^H$
when~\eqref{eqn:1.17.2025.1:14p}  
holds and it places $v \in Y^H$ 
otherwise (with a similar rule applying for when $\bbAbipelem$ last placed $u \in Y^H$).   
Since 
this na\"ive approach may fail to produce a bipartition of $H$, the algorithm $\bbAbipelem$ checks the independence
of $X^H$ and $Y^H$ in $H$, and it is prepared to     
exhaustively  
construct a bipartition $[n] = X_H \dcup Y_H$ of $H$.   
\smallskip 

\begin{framed}  
\noindent  {\sc Algorithm $\bbAbipelem$} 
\begin{enumerate}
\item[$\,$]  {\sc Input:}  an integer $k \geq 3$ and a $k$-graph $H \in \cB_n$.
\smallskip 
\item[$\,$]  {\sc Output:} a bipartition $[n] = X_H \, \dot\cup \, Y_H$ of $H$ in time averaging $O(n^k)$ over $\cB_n$.  
\smallskip 
\item[$\,$]  {\sc Steps:}  
\begin{enumerate}
\item[$\,$]  
\begin{enumerate}
\item[1.]  
Compute in time $O(n^k)$ the $(n-1)$-many joint $H$-degrees 
$\deg_H(i, i+1)$, 
over all 
$1 \leq i \leq n - 1$, where each is computed in time $O(n^{k-1})$.  
\medskip 
\item[2.]  
Put $1 \in X^H$.  Inductively, 
if $i \in X^H$ put 
$\big($cf.~\eqref{eqn:1.17.2025.1:14p}$\big)$  
$$
\qquad 
\qquad 
i+1 \in 
\begin{cases}
X^H  & \text{when $\deg_H(i, i+1) \geq \tfrac{1}{4} \tbinom{n}{k-1} \big(1 - \tfrac{3}{2^k}\big)$,}  \smallskip \\
Y^H  & \text{otherwise,}
\end{cases}
$$
and if $i \in Y^H$ put 
$$
\qquad 
\qquad 
i+1 \in 
\begin{cases}
Y^H  & \text{when $\deg_H(i, i+1) \geq \tfrac{1}{4} \tbinom{n}{k-1} \big(1 - \tfrac{3}{2^k}\big)$,}  \smallskip \\
X^H  & \text{otherwise,}
\end{cases}
$$
which constructs~$[n] = X^H \dcup Y^H$ in time $O(n)$.
\medskip 
\item[3.]   
Check in time $O(n^k)$ whether $X^H$ and $Y^H$ are independent in $H$:    
\medskip 
\begin{enumerate}  
\item[(i)]  when so, set $\big(X_H, Y_H\big) = \big(X^H, Y^H\big)$; 
\medskip  
\item[(ii)] else, exhaustively construct in time $O(2^n n^k)$ a bipartition $[n] = X_H \dcup Y_H$ of $H$.   
\medskip 
\end{enumerate}  
\end{enumerate}  
\end{enumerate}  
\item[$\,$]  {\sc Return:}  $[n] = X_H \dcup Y_H$.  
\end{enumerate}  
\end{framed}  
\smallskip

\subsection*{Proof of Theorem~\ref{thm:main} from $\boldsymbol{\mathbb{A}_{{\rm bip-elem}}}$}  
Fix an integer $k \geq 3$.  
Clearly, {\sc Step}~3 of $\bbAbipelem$ constructs a bipartition $[n] = X_H \dcup Y_H$ of 
every 
$H \in \cB_n$.
It remains to see that it does so in average time $O(n^k)$ over $H \in \cB_n$.  
For that, 
those bipartitions $[n] = X_H \dcup Y_H$ secured by {\sc Step}~3~(i) 
are constructed in maximum time $O(n^k)$.  
We denote by $\cB_{n, {\rm (i)}}$ the set of $H \in \cB_n$ whose bipartition $[n] = X_H \dcup Y_H$
from $\bbAbipelem$ is constructed by {\sc Step}~3~(i).   We prove that 
\begin{equation}
\label{eqn:1.17.2025.3:35p}  
\big| \cB_n \setminus \cB_{n, {\rm (i)}}\big| 
\leq 
|\cB_n| / 2^{\Omega(n^{k-1})}.     
\end{equation}  
If true, the average running time of $\bbAbipelem$ is indeed, 
with $k \geq 3$,  
$$
\tfrac{1}{|\cB_n|} O \Big(
\big|\cB_{n, {\rm (i)}}\big| 
n^k + 
\big| \cB_n \setminus \cB_{n, {\rm (i)}}\big| 2^n n^k
\Big) \leq 
O\Big(n^k \Big(1 + 2^{n - \Omega(n^{k-1})}\Big) \Big) = O(n^k).      
$$

To prove~\eqref{eqn:1.17.2025.3:35p}, we make the following preparations    
for a henceforth arbitrarily fixed 
\begin{equation}
\label{eqn:fixsigma}  
0 < \sigma < 2^{-k}.  
\end{equation}  

\noindent First, we define the following concept based on~\eqref{eqn:fixsigma}.

\begin{definition}[$\boldsymbol{\sigma}$-{\bf standard}]
\label{def:sigmastand}  
\rm 
We define 
$H \in \cB_n$ to be {\it $\sigma$-standard} when all of its bipartitions 
$[n] = \cXH \dcup \cYH$
are {\it $\sigma$-standard w.r.t.~$H$}, meaning that 
all $u \neq v \in [n]$ satisfy 
\begin{equation}
\label{eqn:sigmastand}
\deg_H(u, v) 
= 
\begin{cases}
\tfrac{1}{4} 
\tbinom{n}{k-1} \big(1 - \tfrac{1}{2^{k-1}} \pm \sigma \big)    
& \text{when $u, v \in \cXH$ or $u, v \in \cYH$,} \bigskip \\
\tfrac{1}{4} 
\tbinom{n}{k-1} \big(1 - \tfrac{2}{2^{k-1}} \pm \sigma \big)    
& \text{otherwise.}  
\end{cases}
\end{equation}  
We define 
$\cS_n \subseteq \cB_n$ to be the set of all $\sigma$-standard $H \in \cB_n$.  
\end{definition}

\begin{remark}
In Definition~\ref{def:sigmastand}, it will turn out 
$\big($see the proof of~\eqref{eqn:1.17.2025.3:35p} below$\big)$   
that 
a $\sigma$-standard $H \in \cS_n$ will admit 
a unique (and necessarily $\sigma$-standard) bipartition   
$[n] = \cXH \dcup \cYH$.    
\end{remark}

Second, we will use the following fact.

\begin{fact}
\label{fact:standardchernoff}  
$\big|\cB_n \setminus \cS_n\big| \leq 
|\cB_n| / 2^{\Omega(n^{k-1})}$.    
\end{fact}

\noindent  Fact~\ref{fact:standardchernoff}  
is proved by standard convexity and Chernoff considerations.  For completeness, we sketch this 
proof in the Appendix.    
We now prove~\eqref{eqn:1.17.2025.3:35p}.

\subsubsection*{Proof 
of~\eqref{eqn:1.17.2025.3:35p}}  
By Fact~\ref{fact:standardchernoff},   
it suffices to prove that $\cS_n \subseteq \cB_{n, {\rm (i)}}$.    
For that, fix $H \in \cS_n$, where in what follows we will reference
one of its fixed $\sigma$-standard bipartitions $[n] = \cXH \dcup \cYH$ (but which is 
not 
given to us).   
Recall that 
{\sc Step}~2 of 
$\bbAbipelem$ 
places $1 \in X^H$.  We claim that 
\begin{equation}
\label{eqn:1.17.2025.5:49p}  
\big(X^H, Y^H\big) = 
\begin{cases}
\big(\cXH, \cYH\big) & \text{when $1 \in \cXH$,}  \medskip \\ 
\big(\cYH, \cXH\big) & \text{when $1 \in \cYH$.}  
\end{cases}
\end{equation}  
If true, 
{\sc Step}~3 necessarily confirms the independence of 
$X^H$ and $Y^H$ in $H$
(from that of $\cXH$ and $\cYH$).  Thus, $\bbAbipelem$ 
would 
return $\big(X_H, Y_H\big) = \big(X^H, Y^H\big)$ as the bipartition of $H$, 
which 
places
$H \in \cB_{n, {\rm (i)}}$, as desired.  
In particular, $\bbAbipelem$ recovers the 
arbitrarily 
referenced 
bipartition $[n] = \cXH \dcup \cYH$ of $H \in \cS_n$, which 
makes it the unique bipartition of $H$.

\subsubsection*{Proof 
of~\eqref{eqn:1.17.2025.5:49p}}      
Let, w.l.o.g., 
$1 \in \cXH$, whence $1 \in X^H \cap \cXH$.  
For $1 \leq i \leq n - 1$, 
we prove that $i+1 \in X^H$ if, and only if, $i+1 \in \cXH$.  
Assume, 
w.l.o.g., 
that $i \in \cXH$ 
so that inductively $i \in X^H$.  
If, and on the one hand, $i + 1 \in \cXH$, then 
Definition~\ref{def:sigmastand} guarantees that $i, i+1 \in \cXH$ 
satisfy that 
$$
\deg_H(i, i+1) \geq 
\tfrac{1}{4} \tbinom{n}{k-1}
\big(1 - \tfrac{1}{2^{k-1}} - \sigma\big) 
\stackrel{\eqref{eqn:fixsigma}}{>}    
\tfrac{1}{4} \tbinom{n}{k-1}
\big(1 - \tfrac{3}{2^k} \big),  
$$
so {\sc Step}~2 of $\bbAbipelem$ placed $i+1 \in X^H$.  
If, and on the other hand, $i + 1 \in \cYH = [n] \setminus \cXH$, then 
Definition~\ref{def:sigmastand} guarantees that $i \in \cXH$ and $i+1 \in \cYH$
satisfy that 
$$
\deg_H(i, i+1) \leq 
\tfrac{1}{4} \tbinom{n}{k-1}
\big(1 - \tfrac{2}{2^{k-1}} + \sigma\big) 
\stackrel{\eqref{eqn:fixsigma}}{<}    
\tfrac{1}{4} \tbinom{n}{k-1}
\big(1 - \tfrac{3}{2^k} \big),  
$$
so {\sc Step}~2 of $\bbAbipelem$ placed $i+1 \in Y^H = [n] \setminus X^H$,   
which proves~\eqref{eqn:1.17.2025.5:49p}.

\section{A Second Algorithm for Theorem~\ref{thm:main}}
\label{sec:bbAbipreg}
In this section, we present our regularity-based 
algorithm $\bbAbip$ for 
Theorem~\ref{thm:main}, 
denoted here by $\bbAbipreg$.  
For this, we prepare the following 
constructive 
Szemer\'edi Regularity Lemma~\cite{KRT, Szem1, Szem}.   

\subsection*{Graph regularity}  
Let $G = (V, E)$ be a graph.  For disjoint $\emptyset \neq A, B \subset V$, 
define 
$$
e(A, B)
=   
\big|\big\{ \{a, b\} \in E: \, (a, b) \in A \times B \big\} \big| = d_G(A, B) \cdot |A| |B|, 
$$
where $d_G(A, B)$ is the {\it density} of $(A, B)$.  
For $\eps  > 0$, 
the pair $(A, B)$ is  
{\it $\eps$-regular} 
when 
all $(A', B') \in 2^A \times 2^B$ 
satisfying $|A'| > \eps |A|$ and $|B'| > \eps |B|$ also satisfy $|d_G(A', B') - d_G(A, B)| < \eps$.
A partition 
$\Pi : V = V_1 \, \dot\cup \, \dots \, \dot\cup \, V_t$ is {\it $\eps$-regular} for $G$ when all but 
$\eps t^2$ of its pairs $(V_i, V_j)$, $1 \leq i < j \leq t$, are $\eps$-regular.  
The same partition is {\it $t$-equitable} when $|V_1| \leq \dots \leq |V_t| \leq |V_1| + 1$.    
Szemer\'edi's Regularity Lemma~\cite{Szem1, Szem} guarantees that for every $\eps > 0$, 
every graph $G = (V, E)$
admits an $\eps$-regular and $t$-equitable 
partition 
$\Pi = \Pi_G$  
for $t \leq T_{\rm Szem}(\eps)$, 
where 
$T_{\rm Szem}(\eps)$
is a $O(\eps^{-5})$-level iterated exponential function.    
Alon et al.~\cite{ADLRY} 
established an algorithm which constructs
$\Pi$ 
for $G$ 
in time $O(|V|^{\omega})$, where 
$2 \leq \omega \leq 2.373$ is the exponent of matrix multiplication.  
Kohayakawa, R\"odl, and Thoma~\cite{KRT} optimized the order from~\cite{ADLRY}.

\begin{theorem}[Kohayakawa, R\"odl, and Thoma~\cite{KRT}]  
\label{thm:KRT}
There exists an algorithm $\mathbb{A}_{\rm reg}$ which, 
for all $\eps > 0$, 
constructs of a given graph 
$G = (V, E)$ an $\eps$-regular and $t$-equitable partition $\Pi = \Pi_G$, 
where $t \leq T_{\rm Szem}(\eps)$, in time $O(|V|^2)$.  
\end{theorem}

We also use the following related construct.  
Let a graph $G = (V, E)$ have an 
$\eps$-regular and $t$-equitable partition 
$\Pi:  V = V_1 \, \dot\cup \, \dots \, \dot\cup \, V_t$.  
An 
{\it $(\eps, t)$-cluster} 
$\bC = \bC_{G, \Pi}$ is a graph 
$([t], E_{\bC})$ 
with $e_{\bC} \geq \binom{t}{2} - \eps t^2$ many edges satisfying the rule that, 
for each $1 \leq i < j \leq t$, 
\begin{equation}
\label{eqn:thecluster}
\{i, j\} \in E_{\bC} \qquad \implies \qquad \text{$(V_i, V_j)$ is $\eps$-regular}.  
\end{equation}  

\begin{remark}
\label{rem:12.28.2022}  
Implicit to  Theorem~\ref{thm:KRT}, 
the algorithm $\bbA_{\rm reg}$ also identifies, in time $O(|V|^2)$, 
at least $\tbinom{t}{2} - \eps t^2$ many 
$\eps$-regular pairs $(V_i, V_j)$ of $\Pi = \Pi_G$.  
As such, $\bbA_{\rm reg}$ 
constructs
an $(\eps, t)$-cluster $\bC = \bC_{G, \Pi}$ 
at the same time it constructs 
the $\eps$-regular and $t$-equitable partition $\Pi = \Pi_G$.
\hfill $\Box$  
\end{remark}  

\noindent  Continuing, 
define 
\begin{equation}
\label{eqn:fullpart}
[t]_+ = [t]_{\bC, +} = \big\{i \in [t]: \, \deg_{\bC}(i) \geq \big(1 - 2\eps^{1/2} \big)t \big\}
\qquad \text{and} \qquad 
[t]_- = [t]_{\bC, -} = [t] \setminus [t]_+, 
\end{equation}  
which, in the context of Theorem~\ref{thm:KRT}, 
are identified in time $O(1)$.  
From $e_{\bC} \geq \tbinom{t}{2} - \eps t^2$ follows 
\begin{equation}  
\label{eqn:[t]big}
\big| [t]_+ \big|
\geq \big(1 - 2 \eps^{1/2} \big)t \qquad \text{and} \qquad 
\big|[t]_- \big| \leq 2 \eps^{1/2}t.
\end{equation}  

\subsection*{A sketch of $\boldsymbol{\mathbb{A}_{\text{reg-bip}}}$}  
We first prepare some notation.     
Fix 
a $k$-graph 
$H \in \cB_n$ and a  
$j$-set 
$J \subseteq [n]$ with $j \in \{1, k - 2\}$.    
Define the 
$(k - j)$-uniform 
{\it $J$-link 
hypergraph}  
$H_J = \big([n]\setminus J, E_{H_J}\big)$
by
$E_{H_J} = \big\{ I: \, I \, \dot\cup \, J \in E_H \big\}$. 
When $J = \{v\}$ and $U \subseteq [n] \setminus \{v\}$, we write 
$\big($cf.~\eqref{eqn:basicneigh}$\big)$   
$$
H_v = H_{\{v\}}, 
\quad 
N_H(v) = E_{H_v}, 
\quad 
N_H(v, U) = N_H(v) \cap \tbinom{U}{k-1},  \quad \text{and} \quad 
\deg_H(v, U) = |N_H(v, U)|.
$$
Now, 
for simplicity 
in our sketch, we focus on the case $k = 3$.  
For that, 
select $H \in \cB_n = \cB_n^{(3)}$ uniformly at random.  
Upcoming 
Definition~\ref{def:moderntyp}
and Fact~\ref{fact:oldfact3.2}
will 
detail how 
$H$ admits, with probability 
$1 - \exp \big\{ - \Omega (n^2) \big\}$,   
only 
bipartitions 
$[n] = \cXH \dcup \cYH$
on which it behaves `quasirandomly'.
To hint at
this behavior,  
fix w.l.o.g.~$x \in \cXH$.  
The $x$-link
graph $H_x$ 
will 
inherit from the uniform selection of $H \in \cB_n$ important  
binomial random graph properties  
informally 
described for now very coarsely:   
\smallskip 
\begin{enumerate}
\item[(i)] the induced subgraph $H_x[\cYH]$ `resembles'  
$\mathbb{G}(\cYH, \, 1/2)$; 
\medskip 
\item[(ii)] 
the induced subgraph 
$H_x\big[\cXH \setminus \{x\}\big]
= \mathbb{G}\big(\cXH \setminus \{x\}, \, 0\big)$ is, of course,  
empty; 
\medskip 
\item[(iii)] the induced bipartite subgraph $H_x\big[\cXH \setminus \{x\}, \, \cYH\big]$ `resembles'  
$\mathbb{G}\big(\cXH \setminus \{x\}, \, \cYH; \,  1/2\big)$; 
\smallskip 
\end{enumerate}  
Henceforth 
\begin{multline}  
\label{eqn:outlineassumption}  
\qquad 
\qquad 
\text{{\it we assume that $H_x$ satisfies Properties~(i)--(iii),}} \\
\text{{\it which it will with probability $1 - \exp \big\{ - \Omega(n^2) \big\}$.}}    
\qquad 
\qquad 
\end{multline}

Step~1 of~$\bbAbipreg$
first applies 
Theorem~\ref{thm:KRT} 
to $H_x$ to 
build, in time $O(n^2)$, 
an $\eps$-regular\footnote{ 
For simplicity here, we suppose all pairs $(V_i, V_j)$ above are $\eps$-regular.}  
vertex partition 
$$
V_{H_x} = [n] \setminus \{x\} = 
\big(\cXH \setminus \{x\} \big) \dcup \cYH = 
V_1 \, \dot\cup \, \dots \, \dot\cup \, V_t.      
$$
Its main novelty 
exploits 
that
{\sl each class $V_i$ belongs almost entirely to either $\cXH$ or $\cYH$!}      
In particular, 
\begin{enumerate}
\item[(a)]  
$|V_i \cap \cXH| \leq \eps |V_i|$ or 
$|V_i \cap \cYH| \leq \eps |V_i|$      
lest for some $j$ with $|V_j \cap \cXH| > \eps |V_j|$,  
$$
d_{H_x}\big(V_i \cap \cXH, \, V_j \cap \cXH\big) 
\stackrel{\text{\tiny{(ii)}}}{=} 0 
\qquad \text{and} \qquad 
d_{H_x}\big(V_i \cap \cYH, \, V_j \cap \cXH\big) 
\stackrel{\text{\tiny{(iii)}}}{=}  \tfrac{1}{2} \pm \eps  
$$
contradict the $\eps$-regularity of $(V_i, V_j)$.  
\medskip 
\item[(b)]  
$|V_i \cap \cXH| \leq \eps |V_i|$ is detectable 
in time $O(n^2)$ 
by computing 
$e_{H_x} (V_i)$ and testing  
$$
e_{H_x}(V_i) 
\stackrel{\text{\tiny{(i)}}}{\geq} 
\big(\tfrac{1}{2} - \delta\big) \tbinom{|V_i|}{2} 
\geq 
\tfrac{1}{100} |V_i|^2,  
\qquad \text{for some 
$\delta = \delta (\eps) \to 0^+$ as $\eps \to 0^+$};  
$$
\item[(c)]  
$|V_i \cap \cYH| \leq \eps |V_i|$ is detectable 
in time $O(n^2)$ 
by computing 
$e_{H_x}(V_i)$ and testing  
$$
e_{H_x}(V_i) 
\stackrel{\text{\tiny{(ii)}}}{\leq } 
\gamma \tbinom{|V_i|}{2} 
< 
\tfrac{1}{100} |V_i|^2, \qquad \text{   
for some 
$\gamma = \gamma (\eps) \to 0^+$ as $\eps \to 0^+$}.     
$$
\end{enumerate}  
Step~1 concludes by 
using~(a)--(c)
to 
construct in time $O(n)$ 
a partition 
$[n] = 
\Xapp_H \dcup \Yapp_H$ satisfying  
\begin{equation}
\label{eqn:sketch-close}  
\big|\Xapp_H \, \triangle \, \cXH \big| = 
\big|\Yapp_H \, \triangle \, \cYH \big| \leq 5 \eps^{1/2} n.     
\end{equation}  
Step~2 of $\bbAbipreg$ is then elementary.  
It uses~\eqref{eqn:sketch-close}  
with $\deg_H\big(v, \Xapp_H\big)$ 
and $\deg_H\big(v, \Yapp_H\big)$
over $v \in [n]$ 
to edit $[n] = \Xapp_H \dcup \Yapp_H$ 
in time $O(n^3)$ 
into $[n] = \XedH \dcup \YedH$, which will be   
the unique
bipartition 
$[n] = \cXH \dcup \cYH$
of $H$ satisfying~\eqref{eqn:outlineassumption}.
Since this na\"ive approach may fail
when~\eqref{eqn:outlineassumption} fails, 
$\bbAbipreg$ is prepared to exhaustively construct a bipartition $[n] = X_H \dcup Y_H$ of $H$.  
However, the super-exponential decay in~\eqref{eqn:outlineassumption}  
overcomes 
that exponential construction to preserve 
the average running time $O(n^3)$ in Theorem~\ref{thm:main}.  

It remains to make this outline precise for fixed $k \geq 3$.  
In what follows, 
Step~1 above is essentially 
the subroutine 
$\bbA_{\approx \text{bip}}$    
below 
and 
Step~2 above is 
the subroutine 
$\bbA_{\text{edit}}$ below.    
For the discussion below, 
let $H \in \cB_n$ be fixed arbitrarily, but for the remainder of the paper, 
\begin{multline}
\label{eqn:2.1.2025.3:14p}  
\text{{\sl 
let 
$[n] = \cXH \dcup \cYH$ be an arbitrary but fixed bipartition of $H$}} \\
\text{{\sl which, while not constructed for us, is henceforth assigned to $H$.}}
\end{multline}

\subsection*{The subroutine $\boldsymbol{\bbA}_{\approx \text{\bf bip}}$}  
Informally, 
this subroutine 
intends 
to construct 
in time $O(n^2)$ 
a partition $[n] = \XappH \dcup \YappH$ 
for which
$\big|\XappH \, \triangle \, \cXH\big| = o(n)$ 
or $\big|\XappH \, \triangle \, \cYH| = o(n)$
$\big($cf.~\eqref{eqn:2.1.2025.3:14p}$\big)$.    
Anticipating that 
$\bbA_{\approx {\rm bip}}$ may sometimes fail, 
we refer to $[n] = \XappH \dcup \YappH$ as an {\it approximate bipartition candidate} of $H$.    

\begin{framed}  
\noindent  {\sc Subroutine $\bbA_{\approx \text{{\rm bip}}}$.}
\medskip 
\begin{enumerate}
\item[$\,$]  {\sc Input:}  $(H, J, \eps) \in \cB_n \times \tbinom{[n]}{k-2} \times (0, 1)$. 
\medskip 
\item[$\,$]  {\sc Output:}  an approximate bipartition candidate $[n] = \XappH \dcup \YappH$ of $H$ 
in time $O(n^2)$.
\medskip 
\item[$\,$]  {\sc Steps:}  
\medskip 
\begin{enumerate}
\item[$\,$]  
\begin{enumerate}
\item[1.]  Construct
in time $O(n^2)$ the $J$-link {\sl graph}  
$H_J = (V_{H_J}, E_{H_J})$.  
\medskip 
\item[2.]  Apply Theorem~\ref{thm:KRT} $\big($cf.~Remark~\ref{rem:12.28.2022}$\big)$ to $H_J$ with $\eps$ to   
construct, in time $O(n^2)$, 
\medskip 
\begin{enumerate}
\item[$\circ$]   
an $\eps$-regular and $t$-equitable partition 
$$
\Pi = \Pi_{H_J}: V_{H_J} = V_1 \, \dot\cup \, \dots \, \dot\cup \, V_t, 
$$ 
where $t \leq T_{\rm Szem}(\eps)$; 
\medskip 
\item[$\circ$]  
an $(\eps, t)$-cluster $\bC = \bC_{H_J, \Pi}$.  
\medskip 
\end{enumerate}
\item[3.]  Compute in time $O(n^2)$ 
all $e_{H_J}(V_i)$
over $i \in [t]_+ = [t]_{\bC, +}$ $\big($cf.~\eqref{eqn:fullpart}$\big)$.  
\medskip 
\item[4.]  
Construct in time $O(1)$  
the sets 
\begin{align}  
\qquad 
&[t]_{+, -} = 
[t]_{(\bC, +), \, (J, -)} 
= 
\Big\{ i \in [t]_+: \, e_{H_J}(V_i) < \tfrac{1}{100} |V_i|^2  \Big\} 
\quad \text{and}  \bigskip \nonumber \\
&[t]_{+,+} 
= 
[t]_{(\bC, +), \, (J, +)} 
= \Big\{ j \in [t]_+: \, e_{H_J}(V_j) \geq \tfrac{1}{100} |V_j|^2 \Big\}. 
\label{eqn:t+-++}
\end{align} 
\item[5.]  
Construct in time $O(n)$ 
(cf.~\eqref{eqn:fullpart})  
the sets 
\begin{align}
\qquad 
\qquad 
&\XappH = X^{\approx}_{H, J} = X^{\approx}_{H, J, \eps} = 
J \, \, \dot\cup \, \, 
\Big(\dot\bigcup_{h \in [t]_-} V_h\Big) \, \, 
\dot\cup \, \, 
\Big(
\dot\bigcup_{i \in [t]_{+,-}} V_i\Big)  
\quad 
\text{and} \nonumber \\
&\YappH = Y^{\approx}_{H, J}  
= Y^{\approx}_{H, J, \eps} 
= 
\dot\bigcup_{j \in [t]_{+, +}} V_j.
\label{eqn:9.16.2024.4:41p}  
\end{align}  
\end{enumerate}  
\end{enumerate}  
\item[$\,$]  {\sc Return:}  $[n] = \XappH \dcup \YappH$.  
\end{enumerate}  
\end{framed}  
\smallskip

\noindent We highlight 
details 
from upcoming
Fact~\ref{fact:oldfact3.2}
and 
Lemma~\ref{lem:firstsub}  
for $\bbA_{\approx \text{{\rm bip}}}$.
Fix $(J, \eps) \in \tbinom{[n]}{k-2} \times (0, 1)$ with
$\eps = \eps(k) > 0$ 
from upcoming~\eqref{eqn:Abipconst}.  
We prove 
that, 
for uniformly selected $H \in \cB_n$ $\big($cf.~\eqref{eqn:2.1.2025.3:14p}$\big)$,     
\begin{equation}  
\begin{array}{lll}
\mathbb{P}  \Big[ \, \, \big|\XappH \, \triangle \, \cXH \big| 
= \big|\YappH \, \triangle \, \cYH \big|  
\leq 5 \eps^{1/2} n \, \, \Big| \, \, J \subseteq \cXH \, \, 
\Big]
&  =  & 
1 - 
\exp \big\{ - \Omega(n^2) \big\},  \medskip \\   
\mathbb{P}  \Big[ \, \, \big|\XappH \, \triangle \, \cYH \big| 
 =  \big|\YappH \, \triangle \, \cXH \big|  
\leq 5 \eps^{1/2} n \, \, \Big| \, \, 
J \subseteq \cYH \, \, 
\Big]
& = & 1 - 
\exp \big\{ - \Omega(n^2) \big\},    
\end{array}
\label{eqn:12.27.2022.star}
\end{equation}  
and that, 
for `most' 
$H \in \cB_n$  
$\big($cf.~\eqref{eqn:12.27.2022.star}$\big)$, the 
bipartition 
$[n] = \cXH \, \dot\cup \, \cYH$ 
assigned in~\eqref{eqn:2.1.2025.3:14p}
is unique.
A minor technicality 
with~\eqref{eqn:12.27.2022.star}
when $k \geq 4$
is that, for $H \in \cB_n$ fixed, we can't search for 
$J \in \tbinom{[n]}{k-2}$ satisfying $J \subseteq \cXH$ or $J \subseteq \cYH$ because 
$\cXH \dcup \cYH$ is not given.  
However, we can simply 
iterate $\bbA_{\approx {\rm bip}}$ over all 
$J \in \tbinom{L}{k-2}$
for, w.l.o.g., $L = [2k-5]$, 
since $\big|L \cap \cXH\big| \geq k - 2$ or $\big|L \cap \cYH\big| \geq k - 2$ 
necessarily 
holds.

\subsection*{The subroutine $\boldsymbol{\bbA}_{\rm edit}$}  
Informally, 
this subroutine intends
to edit 
$[n] = \XappH \dcup \YappH$ 
above 
in time $O(n^k)$ 
to arrive at a bipartition $[n] = \XedH \dcup \YedH$ 
of $H$ 
satisfying $\big\{\XedH, \YedH\} = \{\cXH, \cYH\}$  
$\big($cf.~\eqref{eqn:2.1.2025.3:14p}$\big)$.    
Anticipating that 
$\bbA_{\rm edit}$ 
may sometimes fail, we refer to $[n] = \XedH \dcup \YedH$ as a {\it 
bipartition candidate of $H$}.  

\begin{framed}  
\noindent  {\sc Subroutine $\bbA_{\text{{\rm edit}}}$.}
\medskip 
\begin{enumerate}
\item[$\,$]  {\sc Input:}  an $H \in \cB_n$ and a partition $[n] = \XappH \dcup \YappH$
from $\bbA_{\approx {\rm bip}}$.
\medskip 
\item[$\,$]  {\sc Output:}  a 
bipartition candidate  
$[n] = \YedH \dcup \XedH$ of $H$ in time $O(n^k)$.
\medskip 
\item[$\,$]  {\sc Steps:}  
\begin{enumerate}
\item[$\,$]  
\begin{enumerate}
\item[1.]   Compute in time $O(n^k)$ all 
$\deg_H\big(v, \XappH\big)$ and $\deg_H\big(v, \YappH\big)$ 
over $v \in [n]$.  
\medskip 
\item[2.]  
Construct in time $O(n)$ the 
sets 
\begin{equation}
\label{eqn:9.27.2024.12:41p}  
\qquad 
\XedH \, = \, \big\{u \in V: \deg_H\big(u, \XappH\big) \leq \deg_H\big(u, \YappH\big) \big\}
\, = \, [n] \setminus \YedH \medskip 
\end{equation}  
\end{enumerate}
\end{enumerate}  
\item[$\,$]  {\sc Return:}  $[n] = \XedH \dcup \YedH$.  
\end{enumerate} 
\end{framed}  

We highlight details 
from upcoming Fact~\ref{fact:oldfact3.2}
and Lemma~\ref{lem:firstsub}  
for $\bbA_{\rm edit}$.
Fix again 
$\eps = \eps (k)> 0$ 
from upcoming~\eqref{eqn:Abipconst}.  
We prove that,  
for uniformly selected 
$H \in \cB_n$ 
$\big($cf.~\eqref{eqn:2.1.2025.3:14p}  
and~\eqref{eqn:12.27.2022.star}$\big)$,  
\begin{equation}
\label{eqn:12.28.2022.star}
\begin{array}{lll}
\mathbb{P} \Big[ \, \, \big(\XedH, \YedH\big) = \big(\cXH, \cYH\big) \, \, \Big| \, \, 
\big|\XappH \, \triangle \, \cXH\big| = \big|\YappH \, \triangle \, \cYH \big| \leq 
5 \eps^{1/2} n \, \, \Big]
& = & 1 - \exp \big\{- \Omega(n^2) \big\}, \medskip \\  
\mathbb{P} \Big[ \, \, \big(\XedH, \YedH\big) = \big(\cYH, \cXH\big) \, \, \Big| \, \, 
\big|\XappH \, \triangle \, \cYH\big| = \big|\YappH \, \triangle \, \cXH \big| \leq 
5 \eps^{1/2} n \, \, \Big]
& = & 1 - \exp \big\{- \Omega(n^2) \big\}, 
\end{array}
\end{equation}  
and that, 
for `most' $H \in \cB_n$ 
(cf.~\eqref{eqn:12.28.2022.star}), 
the subroutine 
$\bbA_{\rm edit}$ recovers its 
unique bipartition 
$[n] = \cXH \dcup \cYH$  
$\big($whose uniqueness wasn't yet known in~\eqref{eqn:2.1.2025.3:14p}$\big)$.

\subsection*{The algorithm $\boldsymbol{\mathbb{A}_{\text{\bf reg-bip}}}$}  
Henceforth, we fix 
the constant objects 
\begin{equation}
\label{eqn:Abipconst}  
\eps = 
\tfrac{1}{40^2} (4k)^{-2k}, 
\qquad \qquad L = [2k-5], \qquad \qquad  \text{and} 
\qquad \qquad \cJ = \tbinom{L}{k-2}.        
\end{equation}  

\begin{framed}  
\noindent  {\sc Algorithm $\bbA_{\text{{\rm reg-bip}}}$.}
\medskip 
\begin{enumerate}
\item[$\,$]  {\sc Input:}  an $H \in \cB_n$.
\medskip 
\item[$\,$]  {\sc Output:}  a bipartition $[n] = X_H \dcup Y_H$ of $H$ in time averaging $O(n^k)$ over $\cB_n$.  
\medskip 
\item[$\,$]  {\sc Steps:}  
\medskip 
\begin{enumerate}
\item[$\,$]  
\begin{enumerate}
\item[1.]  Apply $\bbA_{\approx {\rm bip}}$:   
\medskip 
\begin{enumerate}
\item[$\circ$]  
apply $\bbA_{\approx\text{bip}}$ 
to $(H, J, \eps)$ over all $J \in \cJ$
$\big($cf.~\eqref{eqn:Abipconst}$\big)$;   
\medskip 
\item[$\circ$]  
construct 
in time $O(n^2)$
the corresponding 
partitions $[n] = X^{\approx}_{H, J} \dcup Y^{\approx}_{H, J}$.  
\medskip 
\end{enumerate}  
\item[2.]  Apply $\bbA_{\rm edit}$:  
\medskip 
\begin{enumerate}
\item[$\circ$]  
apply $\bbA_{\text{edit}}$ to $H$ and 
$[n] = X^{\approx}_{H, J} \dcup Y^{\approx}_{H, J}$  
over all $J \in \cJ$;
\medskip 
\item[$\circ$]  
construct in time $O(n^k)$ the corresponding partitions $[n] = X^{{\rm ed}}_{H, J} \dcup 
Y^{{\rm ed}}_{H, J}$.  
\end{enumerate}  
\medskip 
\item[3.]   
Search in time $O(n^k)$ for $J \in \cJ$ 
s.t.~$X^{{\rm ed}}_{H, J}$ and $Y^{{\rm ed}}_{H, J}$ are independent in $H$:    
\medskip 
\begin{enumerate}  
\item[(i)]  if such $J$ exists, set 
$\big(X_H, Y_H\big) = 
\big(X^{{\rm ed}}_{H, J}, 
Y^{{\rm ed}}_{H, J} \big)$
for any such $J$;
\medskip  
\item[(ii)]  else, exhaustively construct in time $O(2^n n^k)$ a 
bipartition $[n] = X_H \dcup Y_H$ of $H$.   
\medskip 
\end{enumerate}  
\end{enumerate}  
\end{enumerate}  
\item[$\,$]  {\sc Return:}  $[n] = X_H \dcup Y_H$.  
\end{enumerate}  
\end{framed}  

\noindent  This completes the description of the algorithm $\bbAbipreg$.  It remains to use it to prove Theorem~\ref{thm:main}.

\section{A Second Proof of  Theorem~\ref{thm:main}}
\label{sec:mainproofsep}  
In this section, we prove Theorem~\ref{thm:main} 
from the algorithm $\bbAbipreg$ of the previous section.   
Clearly, 
{\sc Step~3} of $\bbAbipreg$ constructs a bipartition $[n] = X_H \dcup Y_H$
of every $H \in \cB_n$.  
It remains to see that it does so in time averaging $O(n^k)$ over $H \in \cB_n$.  
For that, 
those bipartitions $[n] = X_H \dcup Y_H$ secured by its {\sc Step}~3~(i) 
are constructed in maximum time $O(n^k)$.   
We 
denote by $\cB_{n, {\rm (i)}}$ the set of $H \in \cB_n$ 
whose bipartition $[n] = X_H \, \dcup \, Y_H$ from 
$\bbAbipreg$ is constructed by 
{\sc Step}~3~(i).   
We prove that  
\begin{equation}
\label{eqn:ThmAessence}  
\big|\cB_n \setminus \cB_{n, {\rm (i)}}\big| \leq 
|\cB_n| / 2^{\Omega(n^2)}.
\end{equation}  
If true, 
the average running time of $\bbAbip$ is 
indeed 
$$
\tfrac{1}{|\cB_n|}
O\Big(
\big|\cB_{n, {\rm (i)}}\big| n^k 
+ 
\big|\cB_n \setminus \cB_{n, {\rm (i)}}\big| 2^n n^k 
\Big)  
\leq 
O \Big( n^k \Big(1 + 2^{n - \Omega(n^2)} \Big) \Big)
= 
O\big(n^k\big).        
$$

To prove~\eqref{eqn:ThmAessence},   
we make the following preparations for the  
henceforth fixed constant 
\begin{equation}
\label{eqn:thedelta}  
\delta = \eps \big(8 \cdot T_{\rm Szem}(\eps)\big)^{-1},   
\end{equation}      
where 
$\eps > 0$ is from~\eqref{eqn:Abipconst}   
and 
$T_{\rm Szem}(\eps) \in \mathbb{N}$ is guaranteed by the 
Szemer\'edi Regularity Lemma.  
First, we define the following concept
based on~\eqref{eqn:thedelta}.    

\begin{definition}[{\bf $\boldsymbol{\delta}$-typicality}]  
\label{def:moderntyp}  
\rm 
We define  
$H \in \cB_n$ 
to be 
{\it $\delta$-typical} when all of its bipartitions 
$[n] = \cXH \dcup \cYH$ are {\it $\delta$-typical w.r.t.~$H$}, 
meaning they 
satisfy the following properties:
\medskip 
\begin{enumerate}
\item[$\boldsymbol{(P)}$] $\big|\cXH\big|, \big|\cYH\big| = (1\pm \delta)\tfrac{n}{2}$;  
\medskip 
\item[$\boldsymbol{(Q)}$] 
all $(x, y) \in \cXH \times \cYH$ satisfy 
$\deg_H\big(x, \cYH\big) = (1 \pm \delta) \tfrac{1}{2}\tbinom{|\cYH|}{k-1}$ and 
$\deg_H\big(y, \cXH\big) = (1 \pm \delta) \tfrac{1}{2}\tbinom{|\cXH|}{k-1}$; 
\medskip 
\item[$\boldsymbol{(R)}$] 
all pairwise disjoint 
$(J, X_0, Y_0) \in \tbinom{[n]}{k-2}
\times 
2^{\cXH} \times 2^{\cYH}$ 
with 
$|X_0|, |Y_0| > \delta n$ 
satisfy 
\medskip 
\begin{enumerate}
\item[$\boldsymbol{(R_1)}$] 
$e_{H_J}(X_0, Y_0) = (1 \pm \delta) \tfrac{1}{2} |X_0| |Y_0|$;   
\medskip 
\item[$\boldsymbol{(R_2)}$] 
$e_{H_J}(X_0) = (1 \pm \delta) \tfrac{1}{2} \binom{|X_0|}{2}$ when $J \not\subseteq \cXH$ $\big($and $e_{H_J}(X_0) = 0$ otherwise$\big)$; 
\medskip 
\item[$\boldsymbol{(R_3)}$] 
$e_{H_J}(Y_0) = (1 \pm \delta) \tfrac{1}{2} \binom{|Y_0|}{2}$ when $J \not\subseteq \cYH$ $\big($and $e_{H_J}(Y_0) = 0$ otherwise$\big)$.  
\medskip 
\end{enumerate}   
\end{enumerate}  
Define $\cT_n \subseteq \cB_n$ to be the class of 
all $\delta$-typical $H \in \cB_n$.
\end{definition}

\begin{remark}
\rm 
We will see from upcoming Lemma~\ref{lem:firstsub}  
and by the end of the proof of~\eqref{eqn:ThmAessence} (appearing in one moment)  
that every $H \in \cT_n$ admits a unique (and necessarily 
$\delta$-typical) bipartition $[n] = \cXH \dcup \cYH$
which 
$\bbAbipreg$ constructs in maximum time $O(n^k)$.   \hfill $\Box$  
\end{remark}

Second, we will use the following fact.  

\begin{fact}
\label{fact:oldfact3.2}
$|\cB_n \setminus \cT_n| \leq |\cB_n| / 2^{\Omega(n^2)}$. 
\end{fact}

\noindent Fact~\ref{fact:oldfact3.2} is proved by standard convexity and Chernoff considerations.
For completeness, we sketch this proof in the Appendix  
(along with that of the similarly proven 
Fact~\ref{fact:standardchernoff}).       

Third, we fix the assignment, for each 
$(H, J) \in \cT_n \times \cJ$ $\big($cf.~\eqref{eqn:Abipconst}$\big)$,    
\begin{equation}
\label{eqn:assignment}
(H, J) \mapsto 
\big(\, \cXH, \, \, \cYH, \, \,  
X^{\approx}_{H, J}, \, \,  
Y^{\approx}_{H, J}, \, \, 
X^{{\rm ed}}_{H, J}, \, \, 
Y^{{\rm ed}}_{H, J} \, 
\big), 
\end{equation}
where $[n] = \cXH \dcup \cYH$ is the fixed (necessarily $\delta$-typical bipartition) of $H$  
assigned 
in~\eqref{eqn:2.1.2025.3:14p}, 
where 
$[n] = X^{\approx}_{H, J} \dcup  
Y^{\approx}_{H, J}$ is 
from $\bbA_{\approx{\rm bip}}$ for $(H, J, \eps)$
$\big($cf.~\eqref{eqn:Abipconst}$\big)$,   
and where 
$[n] = X^{{\rm ed}}_{H, J} \dcup 
Y^{{\rm ed}}_{H, J}$ 
is then from $\bbA_{\rm edit}$.

Fourth, 
we use the following lemma, 
whose part~(a) is a novelty of the paper.  

\begin{lemma}
\label{lem:firstsub}  
The objects of~\eqref{eqn:assignment}
satisfy 
\begin{enumerate}
\item[(a)]  
$\,$ 
$J \subseteq \cXH$ $\implies$ 
$\big|X^{\approx}_{H, J} \, \triangle \, \cXH\big| 
\leq 5 \eps^{1/2} n$
$\qquad$
and 
$\qquad$
$J \subseteq \cYH$ $\implies$ 
$\big|X^{\approx}_{H, J} \, \triangle \, \cYH\big| 
\leq 5 \eps^{1/2} n$;       
\smallskip 
\item[(b)]  
$\big|X^{\approx}_{H, J} \, \triangle \, \cXH\big| 
\leq 5 \eps^{1/2} n$
$\implies$  
$X^{{\rm ed}}_{H, J} = \cXH$  
$\quad$
and
$\quad$
$\big|X^{\approx}_{H, J} \, \triangle \, \cYH\big| 
\leq 5 \eps^{1/2} n$      
$\implies$  
$X^{{\rm ed}}_{H, J} = \cYH$.   
\end{enumerate}  
\end{lemma}  

\noindent  We defer the proof of 
Lemma~\ref{lem:firstsub}  
to the next section and continue with 
the proof of~\eqref{eqn:ThmAessence}.

\subsection*{Proof of $\boldsymbol{\eqref{eqn:ThmAessence}}$}    
By Fact~\ref{fact:oldfact3.2}, 
it suffices to prove 
that $\cT_n \subseteq \cB_{n, {\rm (i)}}$.    
For that, fix $H \in \cT_n$ and its $\delta$-typical 
bipartition 
$[n] = \cXH \dcup \cYH$ from~\eqref{eqn:assignment}.  
With $L$ and $\cJ$ fixed in~\eqref{eqn:Abipconst}, 
one of $|L \cap \cXH| \geq k-2$ or $|L \cap \cYH| \geq k - 2$ holds, 
so some $J_0 \in \cJ$
satisfies $J_0 \subseteq \cXH$ or $J_0 \subseteq \cYH$.  
Let, w.l.o.g., $J_0 \subseteq \cXH$.  
{\sc Steps~1} and {\sc 2} of 
$\bbAbipreg$ will 
eventually (but 
unwittingly) consider this $J_0 \in \cJ$    
and will construct its corresponding 
partitions $[n] = 
X^{\approx}_{H, J_0} 
\dcup 
Y^{\approx}_{H, J_0}$ and 
$[n]    
= X^{\rm ed}_{H, J_0} 
\dcup 
Y^{\rm ed}_{H, J_0}$   
from~\eqref{eqn:assignment}.  
With 
$J_0 \subseteq \cXH$,   
Lemma~\ref{lem:firstsub}
guarantees 
that 
$\big|X^{\approx}_{H, J_0} \, \triangle \, \cXH \big|  
\leq 5 \eps^{1/2} n$
$\big($by part~(a)$\big)$
and hence 
$X^{{\rm ed}}_{H, J_0} = \cXH$ $\big($by part~(b)$\big)$  
and $Y^{{\rm ed}}_{H, J_0} 
= \cYH$.  Thus,   
{\sc Step~2} 
of $\bbAbip$ 
constructed 
$\big(X^{{\rm ed}}_{H, J_0}, 
Y^{{\rm ed}}_{H, J_0}\big) = (\cXH, \cYH)$ with $J_0$ above.      
Necessarily, 
{\sc Step~3}~(i) 
confirms the independence of 
$X^{{\rm ed}}_{H, J_0} 
= 
\cXH$ and 
$Y^{{\rm ed}}_{H, J_0} 
= 
\cYH$ in $H$, so $\bbAbipreg$
returns $(X_H, Y_H) = 
\big(X^{{\rm ed}}_{H, J_0}, 
Y^{{\rm ed}}_{H, J_0}\big) = (\cXH, \cYH)$, which places    
$H \in \cB_{n, {\rm (i)}}$, 
as desired.  
In particular, 
$\bbAbipreg$ recovers 
the $\delta$-typical bipartition
$[n] = \cXH \dcup \cYH$ of $H \in \cT_n$ 
from~\eqref{eqn:assignment}, which 
makes it 
the unique bipartition of $H \in \cT_n$.

\section{Proof of Lemma~\ref{lem:firstsub}}
\label{sec:firstlemma}  
Fix $(H, J) \in \cT_n \times \cJ$.
With these now fixed, 
we henceforth abbreviate~\eqref{eqn:assignment}  
with 
$$
\cX = \cXH, \quad \cY = \cYH, \quad  
\Xapp = 
X^{\approx}_{H, J}, \quad 
\Yapp = 
Y^{\approx}_{H, J}, \quad   
\Xed = 
X^{\rm ed}_{H, J}, \quad 
\text{and} \qquad 
\Yed = 
Y^{\rm ed}_{H, J}.   
$$
Now, the proof of part~(b) 
of Lemma~\ref{lem:firstsub}
is almost trivial, so we begin our work there.

\subsection*{Proof of part~(b)}  
By symmetry, it suffices to take $|\Xapp \, \triangle \, \cX| \leq 5 \eps^{1/2} n$ 
and prove that $\Xed = \cX$.  To that end, we prove that every fixed $(x, y) \in \cX \times \cY$
satisfies 
\begin{equation}
\label{eqn:9.26.2024.4:59p}
\deg_H\big(x, \Yapp\big) > \deg_H\big(x, \Xapp\big)
\qquad \text{and} \qquad 
\deg_H\big(y, \Xapp\big) > \deg_H\big(y, \Yapp\big),   
\end{equation}
which 
by~\eqref{eqn:9.27.2024.12:41p}  
would place
$(x, y) \in \Xed \times \Yed$ and imply 
$\cX \subseteq \Xed$, $\cY \subseteq \Yed$ 
and 
$(\cX, \cY) = (\Xed, \Yed)$.  
By symmetry, it suffices to argue the first inequality
of~\eqref{eqn:9.26.2024.4:59p}.  
On the one hand, $\deg_H(x, \cX) = 0$ renders 
\begin{equation}
\label{eqn:9.26.2024.5:22p}  
\deg_H\big(x, \Xapp\big) \leq \big|\cY \cap \Xapp\big| \big|\Xapp\big|^{k-2}
\leq 
\big|\Xapp \, \triangle \, \cX\big| n^{k-2} \leq 5 \eps^{1/2} n^{k-1}.  
\end{equation}  
On the other hand, 
$\deg_H\big(x, \Yapp\big) \geq \deg_H\big(x, \, \cY \cap \Yapp\big)$ is at least  
$$
\deg_H (x, \cY) - \big|\cY \cap \Xapp \big| |\cY|^{k-2}  
\geq 
\deg_H (x, \cY) - \big|\Xapp \, \triangle \, \cX\big| n^{k-2}
\geq 
\deg_H (x, \cY) - 5 \eps^{1/2} n^{k-1}.    
$$
Definition~\ref{def:moderntyp} now gives 
\begin{multline*}  
\deg_H\big(x, \Yapp\big) 
+ 5 \eps^{1/2} n^{k-1} 
\geq 
\deg_H\big(x, \cY \big) 
\stackrel{\boldsymbol{(Q)}}{\geq}  
(1 - \delta) \tfrac{1}{2} \tbinom{|\cY|}{k-1} \\
\geq (1 - \delta) \tfrac{1}{2} \Big(\tfrac{|\cY|}{k-1}\Big)^{k-1}   
\stackrel{\boldsymbol{(P)}}{\geq}  
\Big(\tfrac{1 - \delta}{2k} \Big)^k 
n^{k - 1}   
\stackrel{\eqref{eqn:thedelta}}{>} (4k)^{-k} n^{k-1}   
\stackrel{\eqref{eqn:Abipconst}}{=}  
40 \eps^{1/2} n^{k-1}.
\end{multline*}  
Thus, 
$$
\deg_H\big(x, \Yapp\big) 
> 
\big(40\eps^{1/2} - 5 \eps^{1/2} \big) n^{k-1}  
\stackrel{\eqref{eqn:9.26.2024.5:22p}}{>}  
\deg_H\big(x, \Xapp\big).   
$$

\subsection*{Proof of part~(a)}  
By symmetry, 
it suffices to take 
$J \subseteq \cX$ and prove 
$\big|\Xapp \, \triangle \, \cX\big| \leq 5 \eps^{1/2} n$.
We begin 
by discussing some notation.
Let 
$\Xapp$ and $\Yapp$
be constructed
by $\bbAappbip$  
using the 
$\eps$-regular and $t$-equitable partition 
$\Pi = \Pi_{H_J}: V_{H_J} = V_1 \, \dot\cup \, \dots \, \dot\cup \, V_t$ 
and the $(\eps, t)$-cluster graph $\bC = \bC_{H_J, \Pi}$ 
$\big($cf.~\eqref{eqn:thecluster}$\big)$:    
$$
\Xapp
\, \, 
\stackrel{\eqref{eqn:9.16.2024.4:41p}}{=}     
\, \, 
J \, \, \dot\cup \, \, 
\Big(\dot\bigcup_{h \in [t]_-} V_h\Big) \, \, 
\dot\cup \, \, 
\Big(
\dot\bigcup_{i \in [t]_{+,-}} V_i\Big)  
\qquad \text{and} \qquad 
\Yapp
\, \, 
\stackrel{\eqref{eqn:9.16.2024.4:41p}}{=}     
\, \, 
\dot\bigcup_{j \in [t]_{+, +}} V_j        
$$
for $[t] 
= [t]_- \, \dcup \, [t]_+$ and 
$[t]_+ 
= [t]_{+,-} \, \dcup \, [t]_{+,+}$ given by 
$$
[t]_+ 
\stackrel{\eqref{eqn:fullpart}}{=}  
\big\{f \in [t]: \, \deg_{\bC}(f) \geq (1 - 2\eps^{1/2}) t \big\} 
\qquad \text{and} \qquad 
[t]_{+,+} 
\stackrel{\eqref{eqn:t+-++}}{=}  
\Big\{g \in [t]_+: \, e_{H_J}(V_g) \geq \tfrac{1}{100} |V_g|^2
\Big\},  
$$
where 
$\eps > 0$ is from~\eqref{eqn:Abipconst} and   
$t \leq T_{\rm Szem}(\eps)$  
for the constant $T_{\rm Szem}(\eps)$ guaranteed 
by the Szemer\'edi Regularity Lemma.  
For $\ell \in [t]$, we write 
$\cX_{\ell} = \cX \cap V_{\ell}$, $\cY_{\ell} = \cY \cap V_{\ell}$,  
and 
$\big($for the studied
$\Xapp\big)$ 
$$
X^{\approx}_{\ell} 
= \Xapp \cap V_{\ell} 
\stackrel{\eqref{eqn:9.16.2024.4:41p}}{=}     
\left\{
\begin{array}{cc}
\emptyset & \text{if $\ell \in [t]_{+,+}$,}  \\
V_{\ell}  & \text{otherwise.} 
\end{array}
\right.
$$
These considerations establish
that 
\begin{multline*}  
\cX \setminus J 
\, \, 
=
\, \, \bigcup\nolimits_{\ell \in [t]} \cX_{\ell}, 
\qquad 
\Xapp \setminus J 
\, \, 
\stackrel{\eqref{eqn:9.16.2024.4:41p}}{=}     
\, \, 
\bigcup\nolimits_{\ell \in [t]} \Xapp_{\ell}, 
\qquad \text{and}  \\
\Xapp \, \triangle \, \cX \, \, 
\stackrel{\eqref{eqn:9.16.2024.4:41p}}{=}     
\, \,  
\dot\bigcup_{\ell \in [t]} \Big( \big(\Xapp \, \triangle \, \cX\big) \cap V_{\ell} \Big)
\, \, 
= 
\, \, 
\dot\bigcup_{\ell \in [t]} 
\big(\Xapp_{\ell} \, \triangle \, \cX_{\ell}\big),            
\end{multline*}  
where 
every $\ell \in [t]$ satisfies 
$V_{\ell} = \cX_{\ell} \, \dcup \, \cY_{\ell}$ 
and 
$$
X^{\approx}_{\ell} \, \triangle \, \cX_{\ell}  
= 
\left\{
\begin{array}{cc}
\cX_{\ell}  & \text{if $\ell \in [t]_{+,+}$,}  \\
\cY_{\ell}  & \text{otherwise.  } 
\end{array}
\right.
$$
Thus, 
$$
\Xapp \, \triangle \, \cX 
\, \,  =  \, \, 
\Big(
\dot\bigcup_{h \in [t]_-} 
\cY_h 
\Big)
\, \, \dot\cup \, \,       
\Big(
\dot\bigcup_{i \in [t]_{+,-}} 
\cY_i 
\Big) 
\, \, \dot\cup 
\, \, 
\Big(
\dot\bigcup_{j \in [t]_{+,+}} 
\cX_j
\Big)   
$$
satisfies, with 
$|\cY_h| \leq |V_h| \leq \lceil (n - k + 2) / t \rceil
\leq 2 n / t$
for each of 
at most $2 \eps^{1/2} t$
many $h \in 
[t]_-$ 
$\big($cf.~\eqref{eqn:[t]big}$\big)$, 
$$
\big|\Xapp \, \triangle \, \cX \big|
\, \, \leq \, \, 
4 \eps^{1/2} n \, \, + \, \,  
\sum\nolimits_{i \in [t]_{+,-}} 
|
\cY_i 
| 
\, \, + 
\, \, 
\sum\nolimits_{j \in [t]_{+,+}} 
|\cX_j|.  
$$
We claim that each fixed $(i, j) \in [t]_{+,-} \times [t]_{+,+}$ satisfies 
$|\cY_i| < \eps |V_i|$ and $|\cX_j| < \eps |V_j|$, 
which would imply 
$\big| \Xapp \, \triangle \, \cX| \leq \big(4 \eps^{1/2} + \eps\big)n \leq 5 \eps^{1/2}n$.       
More generally, 
we claim that 
each fixed $\ell \in [t]_+$
satisfies 
\begin{enumerate}
\item[(i)]  
$|\cX_{\ell}| < \eps |V_{\ell}|$ or $|\cY_{\ell}| < \eps |V_{\ell}|$,  
\item[(ii)]  
$|\cX_{\ell}| < \eps |V_{\ell}|$ 
$\iff$ $\ell \in [t]_{+,+}$.           
\end{enumerate}    
If true, 
(ii)
gives $|\cX_j| < \eps |V_j|$, as promised.   
Moreover, (ii) also 
gives 
$|\cX_i| \geq \eps |V_i|$, so (i) concludes $|\cY_i| < \eps |V_i|$, as promised.
We now prove~(i) and~(ii).

\subsection*{Proof of~(i)}
Assume, on the contrary, that some fixed 
$\ell \in [t]_+$ satisfies  
\begin{equation}  
\label{eqn:2.7.2022.2:16p}    
|\cX_{\ell}| \geq \eps |V_{\ell}| 
\qquad \text{and} \qquad |\cY_{\ell}| \geq \eps |V_{\ell}|.
\end{equation}  
Clearly, 
\begin{equation}
\label{eqn:someavgj}  
\text{{\sl some 
other 
$\lambda \in N_{\bC}(\ell)$ must also satisfy 
$|\cX_{\lambda}| \geq \eps |V_{\lambda}|$}}
\end{equation}  
lest 
$|J| = k - 2 = O(1)$ and 
$$
|\cX| - |J| 
\, \, = \, \, 
\sum\nolimits_{i \in [t] \setminus N_{\bC}(\ell)} \, \,  |\cX_i| \, \, + \, \, 
\sum\nolimits_{j \in N_{\bC}(\ell)} \, \, |\cX_j|   
\, \, 
\stackrel{\eqref{eqn:fullpart}}{\leq}  
4 \eps^{1/2}n  + 
\eps n 
\stackrel{\eqref{eqn:Abipconst}}{<}   
n / 8 
$$
contradict $|\cX| \geq (1 - \delta) n / 2
\stackrel{\eqref{eqn:thedelta}}{\geq} n / 4$
from 
$\boldsymbol{(P)}$ of 
Definition~\ref{def:moderntyp}.
Now, 
and 
on the one hand,  
$\lambda \in N_{\bC}(\ell)$ 
gives that 
$(V_{\ell}, V_{\lambda})$ 
is $\eps$-regular w.r.t.~$H_J$,  
so~\eqref{eqn:2.7.2022.2:16p}     
and~\eqref{eqn:someavgj} render 
$$
\big| d_{H_J} (\cX_{\ell}, \cX_{\lambda}) - d_{H_J} (V_{\ell}, V_{\lambda}) \big| 
< \eps 
\qquad \text{and} \qquad 
\big| d_{H_J}(\cY_{\ell}, \cX_{\lambda}) - d_{H_J} (V_{\ell}, V_{\lambda}) \big| 
< \eps.   
$$
Crucially, 
\begin{equation}
\label{eqn:9.15.2024.4:43p}
d_{H_J}(\cY_{\ell}, \cX_{\lambda}) 
= 
\big| \underbrace{d_{H_J} (\cX_{\ell}, \cX_{\lambda})}_{\text{ this is 0}}  
- 
d_{H_J}(\cY_{\ell}, \cX_{\lambda}) \big| < 2 \eps   
\stackrel{\eqref{eqn:Abipconst}}{<} 1 / 4
\end{equation}  
because 
\begin{equation}
\label{eqn:9.19.2024.11:26a}
\text{{\sl the bipartition $\cX \, \dot\cup \, \cY$ of $H$ admits no $H$-edges entirely within 
$J \, \dot\cup \, \cX_{\ell} \, \dot\cup \, \cX_{\lambda} \subseteq \cX$.}}  
\end{equation}    
On the other hand, 
$$
|\cY_{\ell}|, |\cX_{\lambda}|
\stackrel{\eqref{eqn:2.7.2022.2:16p}, \,     
\eqref{eqn:someavgj}}{\geq}    
\eps \big\lfloor (n - k + 2) / t \big\rfloor 
\geq 
\eps n / (4t)  
\geq 
\eps \big(4T_{\rm Szem}(\eps)\big)^{-1} n 
\stackrel{\eqref{eqn:thedelta}}{>} \delta n
$$
so 
$\boldsymbol{(R_1)}$ of 
Definition~\ref{def:moderntyp}  
guarantees 
$d_{H_J}(\cY_{\ell}, \cX_{\lambda}) 
\geq (1 - \delta)  / 2 
\stackrel{\eqref{eqn:thedelta}}{\geq}  
1 / 4$,  
which contradicts~\eqref{eqn:9.15.2024.4:43p}.

\subsection*{Proof of~(ii)}  
First, 
let $\ell \in [t]_{+,+}$.  
Then~\eqref{eqn:t+-++} guarantees $e_{H_J}(V_{\ell})  \geq \tfrac{1}{100} |V_{\ell}|^2$  
while~\eqref{eqn:9.19.2024.11:26a}
guarantees 
$e_{H_J}(V_{\ell}) \leq |\cY_{\ell}| |V_{\ell}|$.  Thus, 
$|\cY_{\ell}| \geq \tfrac{1}{100} |V_{\ell}| 
\stackrel{\eqref{eqn:Abipconst}}{\geq} \eps |V_{\ell}|$ so $|\cX_{\ell}| < \eps |V_{\ell}|$
by~(i).   
Second, let 
$|\cX_{\ell}| < \eps |V_{\ell}|$.
Then 
$V_{\ell} = \cX_{\ell} \, \dot\cup \, \cY_{\ell}$ yields
$|\cY_{\ell}| 
> 
(1 - \eps) |V_{\ell}|
\geq n / (8t) 
\stackrel{\eqref{eqn:thedelta}}{\geq} \delta n$,   
so 
$\boldsymbol{(R_3)}$ of 
Definition~\ref{def:moderntyp} guarantees 
$$
e_{H_J}(V_{\ell}) \geq 
e_{H_J}(\cY_{\ell}) 
\geq 
(1 - \delta) \tfrac{1}{2} \tbinom{|\cY_{\ell}|}{2}  
\geq 
(1 - \delta) \tfrac{1}{8}   |\cY_{\ell}|^2 
\geq 
(1 - \delta) (1 - \eps)^2 \tfrac{1}{8}  |V_{\ell}|^2 
\stackrel{\eqref{eqn:Abipconst}, \, \eqref{eqn:thedelta}}{\geq}
\tfrac{1}{100} |V_r|^2,  
$$
which 
places $\ell \in [t]_{+,+}$  
by~\eqref{eqn:t+-++}.   

\section{Concluding Remarks} 
\label{sec:Fano}  
One may use 
the results of this paper and the inequality 
$|\cF_n \setminus \cB_n\big| 
\leq 
\big|\cF_n| / 2^{\Omega(n^2)}$ 
from Theorem~\ref{thm:PSFano} 
to expedite the average running time of its algorithm $\bbA_{\rm Fano}$.
It is perhaps most interesting to do so using the algorithm $\bbAbipelem$
from 
Section~\ref{sec:shortproof}, which we will use
below 
verbatim except for 
appending to 
its exhaustive {\sc Step}~3~(ii) 
(which will fail when $H \in \cF_n \setminus \cB_n$) 
the following suitable enhancement thereof:  
\begin{framed}
$$
\begin{array}{ll}
\text{{\sc Step}~3~(iii)} 
& 
\text{else, exhaustively construct in time $O\big(n^3\sum_{r=1}^n r^n\big) = O(n^{n+4})$ a smallest partition}  \\ 
& 
\text{$[n] = X_{H,1} \dcup \dots \dcup X_{H, r_H}$
into $r_H$ many independent classes of $H$.}
\end{array}
$$
\end{framed}
\noindent  Now, and on the class $\cF_n$, we analyse the performance of $\bbAbipelem$ with the appended 
{\sc Step}~3~(iii).  
For that, 
recall from Definition~\ref{def:sigmastand}   
the class $\cS_n \subseteq \cB_n^{(3)} \subseteq \cF_n$ of $\sigma$-standard $H$, and    
consider now the corresponding partition 
$\cF_n = \cS_n \dcup \big(\cB_n^{(3)} \setminus \cS_n\big) \dcup \big(\cF_n \setminus \cB_n^{(3)}\big)$  
of input $H \in \cF_n$.  

\subsubsection*{Case 1 $(H \in \cS_n)$}  
{\sc Step}~3~(i) of $\bbAbipelem$ constructs in time $O(n^3)$ the unique 
bipartition $[n] = X_H \dcup Y_H$ of $H$ 
$\big($as we proved in~\eqref{eqn:1.17.2025.5:49p}$\big)$.

\subsubsection*{Case 2 $\big(H \in \cB_n^{(3)} \setminus \cS_n\big)$}  
{\sc Step}~3~(ii) constructs in time 
$O(2^n n^3)$ 
a bipartition $[n] = X_H \dcup Y_H$ of $H$ when $E_H \neq \emptyset$, and otherwise realizes
the independence of $[n]$ in $H$ 
in time $O(n^3)$.  

\subsubsection*{Case 3 $\big(H \in \cF_n \setminus \cB_n^{(3)}\big)$}  
{\sc Step}~3~(iii) constructs in time 
$O(n^{n+4})$ a smallest partition $[n] = X_{H,1} \dcup$ $\dots$ $\dcup X_{H, r_H}$ into $r_H \geq 3$ independent classes of $H$.  \\

\noindent  In all three cases, this algorithm 
constructs  
a smallest partition $[n] = X_{H,1} \dcup \dots \dcup X_{H, r_H}$ into $r_H$ independent classes of $H$.  
Its average running time on the class $\cF_n$ is, 
by 
Theorem~\ref{thm:PSFano}  
and Fact~\ref{fact:standardchernoff},   
$$
\tfrac{1}{|\cF_n|}  
O \Big(
n^3 \Big(
|\cS_n| + \big|\cB^{(3)}_n \setminus \cS_n\big| 2^n + 
\big|\cF_n \setminus \cB_n^{(3)}\big| n^{n+1} \Big)
\Big)  
\leq 
O\Big(n^3 \Big(1 + 
2^{n  - \Omega(n^2)}  
+ 
2^{n \log_2 n - \Omega(n^2)}  
\Big)\Big)  
= O(n^3).             
$$

The original steps of $\bbA_{\rm Fano}$ 
in~\cite{PSalm, PSalg} explicitly invoke 
tools from 
hypergraph regularity.  
What is most interesting here is 
that our steps avoid these tools.  The possibility of such an avoidance surprised us.       
In~\cite{LMN}, 
we continue with further work related to~\cite{PSalm, PSalg}.

\section*{Appendix:  Proofs of Facts~\ref{fact:standardchernoff} and~\ref{fact:oldfact3.2}}  
We begin by making several preparations used commonly in the 
proofs of 
Facts~\ref{fact:standardchernoff} and~\ref{fact:oldfact3.2}.    
First, 
fix $\gamma > 0$ and let 
$H \in \cB_n$ be 
{\it $\gamma$-equitable}
when its every bipartition $[n] = X_H \dcup Y_H$ is {\it $\gamma$-equitable}, 
meaning that 
$(1 - \gamma) (n/2) \leq |X_H| \leq |Y_H| \leq (1 + \gamma) (n/2)$.  
Let $\cE_n = \cE_{\gamma, n} \subseteq \cB_n$ 
be the class of all $\gamma$-equitable $H \in \cB_n$.  
For completeness, we 
recall the proof of the following standard result.    

\begin{fact}
\label{fact:convexity}
$|\cB_n \setminus \cE_n| \leq |\cB_n| / 2^{\Omega(n^k)}$.  
\end{fact}  

\begin{proof} 
Fix a partition $[n] = X \, \dot\cup \, Y$ with $x = |X| \leq |Y|$.  
There are\footnote{For $z \in \mathbb{R}$, we use $2 \uparrow z$ to denote $2^z$.}  
precisely   
$2 \uparrow \big( 
\tbinom{n}{k} - \tbinom{x}{k} - \tbinom{n - x}{k} \big)$ 
many $H \in \cB_n$ 
with bipartition 
$[n] = X \, \dot\cup \, Y$.   
  In particular, 
$|\cB_n| \geq 
2 \uparrow \big(\tbinom{n}{k} - \tbinom{\lfloor n/2 \rfloor}{k} - 
\tbinom{\lceil n/2 \rceil}{k} \big)$.  Now,  
$$
\big|\cB_n \setminus \cE_n\big| \leq 
\sum\nolimits_{x = 0}^{(1 - \gamma)n/2} \tbinom{n}{x} 
2 \uparrow \Big(\tbinom{n}{k} - \tbinom{x}{k} - \tbinom{n - x}{k} \Big)
\leq 
\sum\nolimits_{x = 0}^{(1 - \gamma)n/2} 
2 \uparrow \Big(n + \tbinom{n}{k} - \tbinom{x}{k} - \tbinom{n - x}{k} \Big),   
$$
and 
$$
\big|\cB_n \setminus \cE_n\big| \big|\cB_n\big|^{-1}
\leq 
\sum\nolimits_{x = 0}^{(1 - \gamma)n/2} 
2 \uparrow \Big(
n + 
\tbinom{\lfloor n/2 \rfloor}{k} + \tbinom{\lceil n/2 \rceil}{k}
- \tbinom{x}{k} - \tbinom{n - x}{k} \Big).    
$$
Here\footnote{Recall the elementary inequality 
$\binom{x}{k} + \binom{n-x}{k} \geq \binom{x+1}{k} + \binom{n-x-1}{k}$, which holds for all $0 \leq x \leq (n-1)/2$.}, every index $0 \leq x \leq (1 - \gamma)n/2$ 
satisfies 
\begin{align*}  
\tbinom{x}{k} + \tbinom{n - x}{k} - 
\tbinom{\lfloor n/2 \rfloor}{k} - \tbinom{\lceil n/2 \rceil}{k}
- n 
&\geq 
\tbinom{(1 - \gamma) n/2}{k} + \tbinom{(1 + \gamma)n/2}{k} - 
\tbinom{\lfloor n/2 \rfloor}{k} - \tbinom{\lceil n/2 \rceil}{k}
- n \\ 
&= 
\tfrac{(n/2)^k}{k!} \Big( (1 - \gamma)^k + (1 + \gamma)^k - 2 + o(1) \Big)
= \Omega(n^k), 
\end{align*}  
which proves Fact~\ref{fact:convexity}.   
\end{proof}

Second, 
we recall the 
well-known Chernoff inequality~$\big($cf.~\cite{JLR}$\big)$.  

\begin{theorem}
\label{thm:chern}
A binomially distributed random variable $Z \sim \bbB(m, p)$ and $\zeta > 0$ satisfy 
$$
\mathbb{P} \big[ \big| Z - \mathbb{E}[Z] \big| \geq \zeta \mathbb{E}[Z] \big] \leq 
2 \exp \big\{ - \big(\zeta^2 / 3 \big) \mathbb{E}[Z] \big\}.
$$
\end{theorem}

Third, and 
for simplicity of notation,
we write $y = o_x(1)$ for $x, y > 0$ to mean that $y = y(x)$ vanishes with 
$x$ quickly enough to guarantee
some explicitly presented (in context) and elementary inequalities.  

\subsection*{Proof 
of Fact~\ref{fact:standardchernoff}}
Fix $\gamma = o_{\sigma}(1)$ and recall  
$\cE_n = \cE_{\gamma, n}$ 
from Fact~\ref{fact:convexity}.   
We prove that 
\begin{equation}
\label{eqn:1.18.2025.5:58p}  
\big|\cE_n \setminus \cS_n\big| \leq 
|\cB_n| / 2^{\Omega(n^{k-1})},        
\end{equation}
whence 
\begin{equation}  
\label{eqn:1.21.2025.11:09a}  
|\cB_n \setminus \cS_n| \leq |\cB_n \setminus \cE_n| + |\cE_n \setminus \cS_n|
\stackrel{\text{{\tiny 
Fct.\ref{fact:convexity}, \eqref{eqn:1.18.2025.5:58p}}}}{\leq} 
\big(|\cB_n| / 2^{\Omega(n^k)}\big)     
+ 
\big(|\cB_n| / 2^{\Omega(n^{k-1})}\big)     
\leq 
|\cB_n| / 2^{\Omega(n^{k-1})}     
\end{equation}  
proves Fact~\ref{fact:standardchernoff}.

\subsubsection*{Proof of~\eqref{eqn:1.18.2025.5:58p}}  
Fix a $\gamma$-equitable partition $[n] = X \dcup Y$.  In a moment, we will prove that  
\begin{multline*}
\label{eqn:1.22.2025.2:28p}
\text{$[n] = X \dcup Y$ is not $\sigma$-standard 
(cf.~Definition~\ref{def:sigmastand}) w.r.t.~the binomial} \\  
\text{random
$k$-graph $\bbH_{X,Y}$ with probability at most 
$\exp\big\{- \Omega(n^{k-1)}\big\}$.}
\end{multline*}  
Equivalently 
(cf.~Definition~\ref{def:sigmastand}), 
we will prove that 
with probabilty $1 - \exp\big\{- \Omega(n^{k-1)}\big\}$, 
\begin{equation}
\label{eqn:randsigmastand}  
\deg_{\bbH_{X,Y}}(u, v) 
= 
\begin{cases}
\tfrac{1}{4} 
\tbinom{n}{k-1} \big(1 - \tfrac{1}{2^{k-1}} \pm \sigma \big)    
& \text{when $u, v \in X$ or $u, v \in Y$,} \bigskip \\
\tfrac{1}{4} 
\tbinom{n}{k-1} \big(1 - \tfrac{2}{2^{k-1}} \pm \sigma \big)    
& \text{otherwise}  
\end{cases}
\end{equation}  
holds for all $u \neq v \in [n]$.    
Since $\bbH_{X,Y}$ induces the uniform distribution on all subhypergraphs $H$
of the complete bipartite $k$-graph $K^{(k)}[X, Y]$,
we infer that 
$[n] = X \dcup Y$ is $\sigma$-standard w.r.t.~all
but 
$$
2 \uparrow \Big(\tbinom{n}{k} - \tbinom{|X|}{k} - \tbinom{|Y|}{k} - \Omega(n^{k-1}) \Big)
\leq |\cB_n| / 2^{\Omega(n^{k-1})}  
$$
many $H \in \bbH_{X,Y}$.
Now, 
to bound $|\cE_n \setminus \cS_n|$  
we 
simply 
eliminate all $H \in \cE_n$ which admit a $\gamma$-equitable 
but not $\sigma$-standard 
bipartition $[n] = X_H \dcup Y_H$ to infer 
\begin{equation}
\label{eqn:1.22.2025.3:25p}  
|\cE_n \setminus \cS_n| \leq 
2^n |\cB_n| / 2^{\Omega(n^{k-1})} \leq |\cB_n| / 2^{\Omega(n^{k-1})}     
\end{equation}
and hence~\eqref{eqn:1.18.2025.5:58p}.           

To prove~\eqref{eqn:randsigmastand},   
fix $\zeta = o_{\sigma}(1)$ 
$\big($cf.~\eqref{eqn:fixsigma}$\big)$ and $u \neq v \in [n]$.       
Then $\deg_{\bbH_{X,Y}}(u, v)$ is binomially distributed 
with expectation given in~\eqref{eqn:expjointpure}   
and~\eqref{eqn:expjoint} and again below in a moment.    
The Chernoff inequality guarantees 
\begin{multline*}  
\sum\nolimits_{u \neq v \in [n]} 
\bbP \Big[ \Big|
\deg_{\bbH_{X,Y}}(u, v) - 
\bbE\big[\deg_{\bbH_{X,Y}}(u, v)\big]
\Big| 
\geq 
\zeta 
\bbE\big[\deg_{\bbH_{X,Y}}(u, v)\big]
\Big]  \\
\leq 
2
\sum\nolimits_{u \neq v \in [n]}
 \exp \Big\{ - (\zeta^2 / 3) \bbE \big[
\deg_{\bbH_{X,Y}}(u, v)\big] \Big\}
\stackrel{\eqref{eqn:expjoint}}{\leq}    
2n^2 \exp \big\{ - \Omega(n^{k-1})\big\}
\leq \exp \big\{ - \Omega(n^{k-1})\big\}, 
\end{multline*}  
where~\eqref{eqn:expjoint} relies on the $\gamma$-equitability of $[n] = X \dcup Y$.
Thus, with probability $1 - \exp \big\{ - \Omega(n^{k-1})\big\}$, 
\begin{eqnarray*}  
\deg_{\bbH_{X,Y}}(u, v) 
& = &  
(1 \pm \zeta) 
\bbE \big[\deg_{\bbH_{X,Y}}(u, v)\big] \bigskip \\
& \stackrel{\eqref{eqn:expjoint}}{=} &   
(1 \pm \zeta) 
\left\{
\begin{array}{ll}  
\tfrac{1}{4} 
\tbinom{n}{k-1} \big(1 - \tfrac{1}{2^{k-1}} \pm o_{\gamma}(1) \big)    
& 
\text{when $u, v\in X$ or $u, v \in Y$,} \bigskip \\
\tfrac{1}{4} \tbinom{n}{k-1} \big(1 - \tfrac{2}{2^{k-1}} \pm o_{\gamma}(1) \big)    
& \text{otherwise}
\end{array}
\right.
\end{eqnarray*}  
holds for all $u \neq v \in [n]$, 
which yields~\eqref{eqn:randsigmastand}  
with $\gamma, \zeta = o_{\sigma}(1)$ sufficiently small.

\subsection*{Proof of Fact~\ref{fact:oldfact3.2}}  
Respectively 
in Definition~\ref{def:moderntyp},     
let $\cP_n = \cP_{\delta, n}$, 
$\cQ_n = \cQ_{\delta, n}$, 
and 
$\cR_n = \cR_{\delta, n}$ 
be the class of all $H \in \cB_n$ whose every bipartition $[n] = X_H \dcup Y_H$ satisfies, respectively, 
Properties~$\boldsymbol{(P)}$,  
$\boldsymbol{(Q)}$, and~$\boldsymbol{(R)}$.     
Then $\cP_n = \cE_{\delta, n}$ is the class of all $\delta$-equitable $H \in \cB_n$, 
where Fact~\ref{fact:convexity} 
(which is based on convexity)  
guarantees $|\cB_n \setminus \cP_n| = 
|\cB_n \setminus \cE_{\delta, n}| \leq |\cB_n| / 2^{\Omega(n^k)}$.  
We will prove  
\begin{eqnarray}
\big|\cB_n \setminus \cQ_n \big| & \leq & |\cB_n| / 2^{\Omega(n^{k-1})}, \label{eqn:1.21.2025.11:34a} \bigskip \\
\big|\cB_n \setminus \cR_n \big| & \leq & |\cB_n| / 2^{\Omega(n^2)}, \label{eqn:1.21.2025.11:35a} 
\end{eqnarray}  
whence 
\begin{multline*}  
\big|\cB_n \setminus \big(\cP_n \, \cap \, \cQ_n \, \cap \, \cR_n \big) \big|
\leq 
\big|\cB_n \setminus \cP_n\big|
+
\big|\cB_n \setminus \cQ_n\big|
+
\big|\cB_n \setminus \cR_n\big|  \\
\leq 
\big(|\cB_n| / 2^{\Omega(n^k)}\big)
+   
\big(|\cB_n| / 2^{\Omega(n^{k-1})}\big)
+   
\big(|\cB_n| / 2^{\Omega(n^2)}\big)
\leq 
|\cB_n| / 2^{\Omega(n^2)}
\end{multline*}   
holds 
with $k \geq 3$, 
proving Fact~\ref{fact:oldfact3.2}.

\subsubsection*{Proof of~\eqref{eqn:1.21.2025.11:34a}}
With $\cP_n = \cE_{\delta, n}$
considered in Fact.\ref{fact:convexity}, 
we prove that 
\begin{equation}
\label{eqn:new1.21.2025.11:34a} 
\big|\cP_n \setminus \cQ_n \big| 
\leq |\cB_n| / 2^{\Omega(n^{k-1})},    
\end{equation} 
whence 
$$
|\cB_n \setminus \cQ_n| \leq |\cB_n \setminus \cP_n| + |\cP_n \setminus \cQ_n|
\stackrel{\text{\tiny{Fct.\ref{fact:convexity}, 
\eqref{eqn:new1.21.2025.11:34a}}}}{\leq}  
\big(|\cB_n| / 2^{\Omega(n^k)}\big)
+ 
\big(|\cB_n| / 2^{\Omega(n^{k-1})}\big)
\leq 
|\cB_n| / 2^{\Omega(n^{k-1})}
$$
proves~\eqref{eqn:1.21.2025.11:34a}.  
Now, the proof 
of~\eqref{eqn:new1.21.2025.11:34a}
is 
analogous 
to that of~\eqref{eqn:1.18.2025.5:58p}.
It suffices $\big($cf.~\eqref{eqn:1.22.2025.3:25p}$\big)$  
to fix a $\delta$-equitable partition $[n] = X \dcup Y$
and prove that 
$$
\text{$\bbH_{X,Y}$ satisfies 
Property~$\boldsymbol{(Q)}$   
of Definition~\ref{def:moderntyp} with probability $1 - \exp \big\{- \Omega(n^{k-1})\big\}$}.  
$$
Equivalently
(cf.~Definition~\ref{def:moderntyp}), 
we prove that 
with probability $1 - \exp \big\{- \Omega(n^{k-1})\big\}$, 
\begin{equation}
\label{eqn:1.22.2025.11:18a}  
\deg_{\bbH_{X,Y}}(x, Y) = (1 \pm \delta) 
\tfrac{1}{2} \tbinom{|Y|}{k-1}
\qquad \text{and} \qquad 
\deg_{\bbH_{X,Y}}(y, X) = (1 \pm \delta) 
\tfrac{1}{2} \tbinom{|X|}{k-1}
\end{equation}
hold 
for all $(x, y) \in X \times Y$. 
But~\eqref{eqn:1.22.2025.11:18a}  
follows immediately from the Chernoff inequality.  Indeed, 
such 
$\deg_{\mathbb{H}_{X,Y}}(x, Y)$
and 
$\deg_{\mathbb{H}_{X,Y}}(y, X)$ are binomially distributed 
with respective expectations 
$$
\bbE \big[\deg_{\mathbb{H}_{X,Y}}(x, Y) \big]
= 
\tfrac{1}{2} \tbinom{|Y|}{k-1} = \Omega(n^{k-1}) 
\qquad \text{and} \qquad 
\bbE \big[\deg_{\mathbb{H}_{X,Y}}(y, X) \big]
= 
\tfrac{1}{2} \tbinom{|X|}{k-1} = \Omega(n^{k-1}).      
$$
Thus,  
the Chernoff inequality guarantees 
with $\{A, B\} = \{X, Y\}$ that 
$$
\sum\nolimits_{a \in A}
\mathbb{P} \Big[  \, \Big|  
\deg_{\mathbb{H}_{X,Y}}(a, B) - 
\tfrac{1}{2} \tbinom{|B|}{k-1} 
\Big|
\geq 
\tfrac{\delta}{2} \tbinom{|B|}{k-1} 
\Big]
\leq 
2n \exp \Big\{ - \big(\delta^2 / 6\big) \tbinom{|B|}{k-1} \Big\}
\leq 
\exp \big\{ - \Omega(n^{k-1}) \big\}    
$$
from which~\eqref{eqn:1.22.2025.11:18a} follows.   

\subsubsection*{Proof of~\eqref{eqn:1.21.2025.11:35a}}   
Similarly 
to~\eqref{eqn:1.18.2025.5:58p} 
and~\eqref{eqn:1.21.2025.11:34a},
it suffices to 
fix a partition $[n] = X \dcup Y$ 
and prove that 
$$
\text{$\bbH_{X,Y}$ satisfies 
Property~$\boldsymbol{(R)}$   
of Definition~\ref{def:moderntyp} with probability $1 - \exp\{- \Omega(n^2)\}$.}
$$
Equivalently $\big($cf.~Definition~\ref{def:moderntyp}$\big)$, 
we prove that for each $i = 1, 2, 3$, 
\begin{equation}
\label{eqn:1.22.2025.1:51p}
\text{$\bbH_{X,Y}$ satisfies 
Property~$\boldsymbol{(R_i)}$   
of Definition~\ref{def:moderntyp} with probability $1 - \exp\{- \Omega(n^2)\}$,}
\end{equation}
where these properties are slighly clumsy.  
For~\eqref{eqn:1.22.2025.1:51p}, 
fix a pairwise disjoint triple $(J, X_0, Y_0) \in \tbinom{[n]}{k-2} \times 2^X \times 2^Y$ satisfying $|X_0|, |Y_0| > \delta n$.    
Each of 
$e_{\bbH^{X,Y}_J}(X_0, Y_0)$, 
$e_{\bbH^{X,Y}_J}(X_0)$, 
and 
$e_{\bbH^{X,Y}_J}(Y_0)$ 
is binomially distributed with respective expectations 
$$
\begin{array}{lll}  
\bbE\big[
e_{\bbH^{X,Y}_J}(X_0, Y_0)\big]
& = &  
\tfrac{1}{2} |X_0| |Y_0|, \bigskip \\
\bbE\big[
e_{\bbH^{X,Y}_J}(X_0)\big]
& = &  
\begin{cases}
\tfrac{1}{2} \tbinom{|X_0|}{2}               & \text{when $J \not\subseteq X$,}  \bigskip \\
0                                            & \text{otherwise,} 
\end{cases} \bigskip \\
\bbE \big[
e_{\bbH^{X,Y}_J}(Y_0)\big]
& = &  
\begin{cases}
\tfrac{1}{2} \tbinom{|Y_0|}{2}               & \text{when $J \not\subseteq Y$,}  \bigskip \\
0                                            & \text{otherwise.} 
\end{cases}  
\end{array}
$$
To prove~\eqref{eqn:1.22.2025.1:51p}
with $i = 1$, 
we sum over all 
pairwise disjoint triples $(J, X_0, Y_0) \in \tbinom{[n]}{k-2} \times 2^X \times 2^Y$ satisfying $|X_0|, |Y_0| > \delta n$    
and use 
the Chernoff inequality to infer 
that 
\begin{align*}  
\sum
\bbP\Big[
\big|
e_{\bbH^{X,Y}_J}(X_0, Y_0)
- 
\tfrac{1}{2} |X_0| |Y_0| 
\big|
\geq \tfrac{\delta}{2} |X_0| |Y_0|\Big]  
&\leq 
2 
\sum
\exp \big\{ - (\delta^2 / 6) |X_0| |Y_0| \big\}  \\
&\leq 2 n^{k-2} 4^n 
\exp \big\{ - (\delta^4 / 6) n^2 \big\}
\leq 
\exp \big\{ - \Omega(n^2) \big\}.    
\end{align*}  
To prove~\eqref{eqn:1.22.2025.1:51p}
with $i = 2$, 
we sum over all 
pairwise disjoint triples $(J, X_0, Y_0) \in \tbinom{[n]}{k-2} \times 2^X \times 2^Y$ satisfying 
$J \not\subseteq X$ and 
$|X_0|, |Y_0| > \delta n$ and use  
the Chernoff inequality to infer 
that 
\begin{align*}  
\sum
\bbP\Big[\,
\Big|
e_{\bbH^{X,Y}_J}(X_0)
- 
\tfrac{1}{2} \tbinom{|X_0|}{2}
\Big|
\geq \tfrac{\delta}{2} \tbinom{|X_0|}{2}\Big]
&\leq 
2 
\sum
\exp \Big\{ - (\delta^2 / 6) \tbinom{|X_0|}{2} \Big\}  \\
&\leq 2 n^{k-2} 4^n 
\exp \big\{ - \big(\delta^4 / (24)\big) n^2 \big\}
\leq 
\exp \big\{ - \Omega(n^2) \big\}.     
\end{align*}  
To prove~\eqref{eqn:1.22.2025.1:51p}
with $i = 3$, 
we 
sum over all 
pairwise disjoint triples $(J, X_0, Y_0) \in \tbinom{[n]}{k-2} \times 2^X \times 2^Y$ satisfying 
$J \not\subseteq Y$ and 
$|X_0|, |Y_0| > \delta n$ and use   
the Chernoff inequality to infer 
that 
\begin{align*}  
\sum
\bbP\Big[\,
\Big|
e_{\bbH^{X,Y}_J}(Y_0)
- 
\tfrac{1}{2} \tbinom{|Y_0|}{2}
\Big|
\geq \tfrac{\delta}{2} \tbinom{|Y_0|}{2}\Big]
&\leq 
2 
\sum 
\exp \Big\{ - (\delta^2 / 6) \tbinom{|Y_0|}{2} \Big\}  \\
&\leq 2 n^{k-2} 4^n 
\exp \big\{ - \big(\delta^4 / (24)\big) n^2 \big\}
\leq 
\exp \big\{ - \Omega(n^2) \big\}.     
\end{align*}

\end{document}